\documentclass[a4paper,11pt]{article}
\usepackage{latexsym}

\newcommand{\thm}{\sc{Theorem}\ }
\newcommand{\pro}{\sc{Proposition}\ }
\newcommand{\cor}{\sc{Corollary}\ }
\newcommand{\exa}{\sc{Example}\ }
\begin{document}
\begin{center}
  {\bf \Large Curvilinear coordinates on generic conformally flat hypersurfaces and
constant curvature 2-metrics}\\
\vspace*{4mm}
	Francis E.~BURSTALL,  Udo HERTRICH-JEROMIN  and  Yoshihiko SUYAMA
\end{center}

{\bf Abstract} \ \ There is a one-to-one correspondence between
associated families of generic conformally flat (local-)hypersurfaces
in 4-dimensional space forms and conformally flat 3-metrics with the
Guichard condition. In this paper, we study the space of conformally
flat 3-metrics with the Guichard condition: \ for a conformally flat
3-metric with the Guichard condition in the interior of the space, an
evolution of orthogonal (local-)Riemannian $2$-metrics with constant Gauss curvature $-1$ is determined; 
for a $2$-metric belonging to a certain class of orthogonal analytic
$2$-metrics with constant Gauss curvature $-1$, a one-parameter family
of conformally flat 3-metrics with the Guichard condition is
determined as evolutions issuing from the $2$-metric.

\footnote[0]{The first and second authors were partly supported by the Bath
Institute for Mathematical Innovation (July 2015).}
\footnote[0]{The second author was partly supported by JSPS-FWF Joint Research Project (the period: 2014-2016).}
\footnote[0]{The third author was partly supported by Grant-in-Aid for Scientific Research (C) (26400076) and JSPS-FWF Joint Research Project (the period: 2014-2015).}

\section*{Introduction}
\label{sec:introduction}

The aim of this paper is to study the space of generic conformally flat (local-)hypersurfaces of dimension 3 in 4-dimensional space forms via conformally flat 3-metrics with the Guichard condition. 
Here, a hypersurface is called generic if it has distinct
principal curvatures at each point.

\vspace*{2mm} A complete local classification of conformally
flat hypersurfaces in n-dimensional space forms, $n\geq 5$,
was given by Cartan\cite{ca}: a hypersurface in an
$n$-dimensional space form, $n\geq 5$, is conformally flat if
and only if it is a branched channel hypersurface, i.e., if
and only if it is quasi-umbilic.  3-dimensional branched
channel hypersurfaces in a 4-dimensional space form are
known to be conformally flat as well, but there are also
generic 3-dimensional conformally flat hypersurfaces.  To
find the complete (local) classification of these
hypersurfaces is an open problem.  However, several partial
classification results of generic conformally flat
hypersurfaces were given in \cite{hs1}, \cite{hs2},
\cite{su1} (and see also \cite{su2} and \cite{su3}).  In
this paper, we relate generic conformally flat hypersurfaces
to families of orthogonal (local-)Riemannian $2$-metrics with
constant Gauss curvature $-1$.

Any generic conformally flat hypersurface in a 4-dimensional space form has a special curvilinear coordinate system $(x,y,z)$ satisfying the following conditions: 

(1) all coordinate lines are principal curvature lines. 

(2) its first fundamental form $I$ is expressed as 
$$I=l_1^2(dx)^2+l_2^2(dy)^2+l_3^2(dz)^2.$$

(3) the functions $l_i^2$ ($i=1,2,3$) satisfy a Guichard
condition $l_i^2+l_j^2=l_k^2$, where $\{i,j,k\}$ is some \hspace*{0.5cm} permutation of $\{1, 2, 3\}$. 

Such a coordinate system is called a principal Guichard net of a generic conformally flat hypersurface. 
We note that the Guichard condition ((2) and (3)) is conformally invariant, that is, it is preserved under conformal changes of the induced metric. 
Therefore, a principal Guichard net of a generic conformally flat hypersurface in a 4-dimensional space form can be mapped to Euclidean 3-space $R^3$ using a conformal coordinate system of the hypersurface to obtain a Guichard net in $R^3$, which is unique up to M\"{o}bius transformation. 
Thus, we can recognise that a Guichard net is a pair \ $\{(x,y,z), [g]\}$ \ of a coordinate system $(x,y,z)$ on a simply connected domain $U$ in $R^3$ and the conformal class $[g]$ of a conformally flat metric $g$ satisfying the Guichard condition with respect to the coordinate system.   

Conversely, for a given Guichard net \ $\{(x^1,x^2,x^3), [g]\}$, \  there exists 
a generic conformally flat hypersurface with its canonical principal
Guichard net in a 4-dimensional space form, uniquely up to M\"{o}bius
transformation (cf. \cite{he3} \S2.4.6). 
Here, the term ``canonical Guichard net'' refers to the conditions \
$\theta^1=dx$, \ $\theta^2=dy$ \ and \ $\theta^3=dz$ \ for the
conformal fundamental 1-forms $\theta^i$  ($i=1,2,3$) of the
hypersurface (cf. \cite{he3} \S2.3.3). 
Then, the coordinates $x,y,z$ are determined up to sign and constant of integration, as $\theta^i$ ($i=1,2,3$) are only determined up to sign. Here, we assume that the domain $U$, where $g$ is defined, intersects the plane $z=0$ for the sake of simplicity for the description later. 
This existence theorem was obtained by study of the integrability condition on a generic conformally flat hypersurface with the canonical principal Guichard net in the conformal 4-sphere. 
A method to determine the the first and the second fundamental forms for a generic conformally flat hypersurface realised in $R^4$ from a Guichard net has been provided in \cite{hs3}.

Certain non-trivial transformations (resp. deformations) act on the
space of generic conformally flat hypersurfaces: each hypersurface has
an associated family, which is a one-parameter family of
non-equivalent generic conformally flat hypersurfaces with the same
Guichard net (cf. \cite{fp}, see also \cite{hs2} and
\cite{su2}, or \cite{bc3} for a more general statement); 
each hypersurface in $R^4$ has its dual generic conformally flat hypersurface in $R^4$, which generally belongs to a different conformal class (or has a different Guichard net) from the one of the original hypersurface (cf. \cite{hsuy}, \cite{bc2}), but, as to its principal coordinate system determined from the Guichard net, we can take the same coordinate system as in the original hypersurface (cf. \cite{hsuy}). 

Let $\iota_p$ be an inversion acting on $R^4$ with respect
to 3-sphere $S^3_p$ of radius $1$ and center $p$.  For a
generic conformally flat hypersurface $f$ in $R^4$, both
duals $(\iota_p f)^*$ and $(\iota_q f)^*$ of $\iota_p f$ and
$\iota_q f$, respectively, are generally non-equivalent if
$p\neq q$ (cf. \cite{hsuy}).  Hence, a five dimensional set
of generic conformally flat hypersurfaces is constructed
from one hypersurface (see \cite{bc2} for another proof of
this fact). When we further consider
$(\iota_q(\iota_p f)^*)^*$ and so on, the space of generic
conformally flat hypersurfaces seems to be very large.

Let $\kappa_i$ ($i=1,2,3$) be the principal curvatures corresponding to the coordinate lines $x$, $y$ and $z$, respectively, of a generic conformally flat hypersurface, and for the sake of simplicity suppose that $\kappa_3$ is the middle principal curvature for the hypersurface, i.e., \ $\kappa_1>\kappa_3>\kappa_2$ \ or \ $\kappa_1<\kappa_3<\kappa_2$. \ 
Then, by the Guichard condition there is a function $\varphi=\varphi(x,y,z)$ such that a metric $g$,
$$g=\cos^2\varphi(dx)^2+\sin^2\varphi(dy)^2+(dz)^2, \eqno{(1)}$$
together with the coordinate system $(x,y,z)$ is a
representative of the Guichard net determined by the
hypersurface.

Thus, the existence problem of generic conformally flat
hypersurfaces is reduced to that of conformally flat metrics
$g$ (resp. functions $\varphi$) given by (1).

\vspace{2mm}
Now, we assume that all metrics $g$ given by (1) (resp. all hypersurfaces) are of $C^{\infty}$-class.
Let $\varphi_{z}$ (resp. $\varphi_{xz}$) be the first derivative (resp. the second derivative) of $\varphi$ with respect to $z$ (resp. with respect to $x$ and $z$). Our main Theorem 1 is as follows (see Theorem 1 in \S1 and Theorem 2 in \S2.1):\\

{\sc Main} {\thm}1. {\it \ \ Let $g$ be a conformally flat 3-metric defined by (1) from a function $\varphi(x,y,z)$. Then, we have the following facts (1) and (2): 

(1) \ There is a function $\psi(x,y,z)$ such that \ 
$\psi_{xz}=-\varphi_{xz}\cot\varphi, \ \ \psi_{yz}=\varphi_{yz}\tan\varphi.$

(2) \ Suppose that \ $\varphi_{xz}\neq 0$ \ and \ $\varphi_{yz}\neq 0$ are satisfied. Let us define functions ${\hat A}(x,y,z)$ and ${\hat B}(x,y,z)$ by
$${\hat A}:=-\displaystyle\frac{\varphi_{xz}}{\varphi_z\sin\varphi}=\displaystyle\frac{\psi_{xz}}{\varphi_z\cos\varphi}, \hspace{1cm} {\hat B}:=\displaystyle\frac{\varphi_{yz}}{\varphi_z\cos\varphi}=\displaystyle\frac{\psi_{yz}}{\varphi_z\sin\varphi}.$$  
Then, the Riemannian $2$-metric \ ${\hat g}(z):={\hat A}^2(x,y,z)(dx)^2+{\hat B}^2(x,y,z)(dy)^2$ \ for any $z$ has constant Gauss curvature $K_{{\hat g}(z)}\equiv -1$.}\\

When $\varphi$ in a conformally flat 3-metric $g$ satisfies the conditions \ $\varphi_{xz}=\varphi_{yz}=0$, \ $g$ leads to a generic conformally flat hypersurface either of product-type or with cyclic Guichard net. 
For hypersurfaces of product-type, see  (\cite{su2}, \S2.2) and \cite{la}. 
All generic conformally flat hypersurfaces with cyclic Guichard net were explicitly realised in 4-dimensional space forms and completely classified in \cite{hs1}. By the Main Theorem 1, we know that two kinds of hypersurfaces of product-type and with cyclic Guichard net determined from $\varphi$ satisfying $\varphi_{xz}=0$ and $\varphi_{yz}=0$ lie in the boundary of the space of generic conformally flat hypersurfaces.  

\vspace{2mm}
Next, let \ ${\hat g}={\hat A}^2(x,y)(dx)^2+{\hat
  B}^2(x,y)(dy)^2$ \ be a Riemannian $2$-metric with constant
Gauss curvature $-1$ defined on a simply connected domain
$V$ in the $(x,y)$-plane.
Then, there are three functions $\varphi(x,y)$, $\varphi_z(x,y)$ and $\psi_z(x,y)$ on $V$ satisfying the following condition:
$${\hat A}=-\displaystyle\frac{\varphi_{zx}}{\varphi_z\sin\varphi}=\displaystyle\frac{\psi_{zx}}{\varphi_z\cos\varphi}, \hspace{1cm}{\hat B}=\displaystyle\frac{\varphi_{zy}}{\varphi_z\cos\varphi}=\displaystyle\frac{\psi_{zy}}{\varphi_z\sin\varphi}.$$
In these equations, $\varphi(x,y)$ is uniquely determined from ${\hat g}$ by giving $\varphi(0,0)=\lambda$, but $\varphi_z(x,y)$ and $\psi_z(x,y)$ are only determined up to the same constant multiple $c\neq 0$ even if we assume $\psi_z(0,0)=0$, that is, \ $\varphi_z(x,y)=\varphi^c_z(x,y):=c\varphi^1_z(x,y)$ \ and \ $\psi_z(x,y)=\psi^c_z(x,y):=c\psi^1_z(x,y)$ (see Theorem 3 in \S2.2). \ 

In \S4, we study the following system of evolution equations in $z$,
$$
\begin{array}{l}
\psi_{zz}=(\varphi_{xx}-\varphi_{yy})\sin2\varphi-(\psi_{xx}-\psi_{yy})\cos2\varphi, \\
\varphi_{zz}=(\varphi_{xx}-\varphi_{yy})\cos2\varphi+(\psi_{xx}-\psi_{yy})\sin2\varphi. \\
\end{array}
\eqno{(2)}
$$
In \S1, Theorem~1, we show that the functions $\varphi$,
$\psi$ arising from a Guichard net as in Main Theorem 1 are
solutions of the system (2) and investigate whether the
converse is true.  The Cauchy--Kovalevskaya theorem ensures
that solutions of (2) exist for given real-analytic initial data
$\varphi(x,y)$, $\varphi_z(x,y)$, $\psi(x,y)$ and
$\psi_z(x,y)$ on the coordinate surface $z=0$.  As we have seen,
this data gives rise to a constant Gauss curvature metric
$\hat{g}$ but additional equations are required on that data
for the corresponding solution of (2) to give rise to a
Guichard net (see \S4, Proposition~4.2) and so an evolution
$\hat{g}(z)$ of constant curvature $2$-metrics.  In
particular, not all such $\hat{g}$ can serve as the initial
metric for such an evolution (see Example~2 in \S3.2).

In general, the necessary equations on initial data are
complicated and difficult to understand (see
Proposition~3.2).  However, some simplification can be
achieved by requiring that these equations are satisfied for
\emph{all} initial data giving rise to the same $2$-metric
$\hat{g}$, that is, for $\varphi(x,y)$, $\varphi^c_z(x,y)$
and $\psi^c_z(x,y)$, for all $c\neq 0$.  In this situation,
we can describe the requirements on initial conditions to
get an evolution on $2$-metrics and then a $1$-parameter family
$g^c$ of $3$-metrics providing Guichard nets.  This is the
content of Main Theorem 2 which we now state.

\vspace{2mm}
Let \ $Lf=(Lf)(x,y)=(f_{xx}-f_{yy})(x,y)$ \ for a function $f=f(x,y)$ and \ $\varphi_z(x,y)=\varphi^c_z(x,y):=c\varphi^1_z(x,y)$. \ Our main Theorem 2 is as follows (see Theorems 5, 6 in \S3.2 and Theorem 7 in \S4). \\

{\sc Main} {\thm}2. \ \ {\it Let two classes (A) and (B)
  of pairs of functions $\varphi(x,y)$ and
  $\varphi^1_z(x,y)$ be defined as follows:

(A) $\varphi(x,y)$ and $\varphi^1_z(x,y)$ are given by 
$$\cos^2\varphi(x,y):=\frac{1}{1+e^{D(y)}}, \hspace{1cm} (\varphi^1_z)^2(x,y):=\zeta(x)\sin^2\varphi(x,y),$$
respectively, with non-constant {\rm analytic} functions
$\zeta(x)$, $D(y)$ of one-variable. Similarly, functions
$\varphi(x,y)$ determined by \
$\cos^2\varphi(x,y):=1/(1+e^{C(x)})$ \ are also included in
this class, then the partners $\varphi^1_z(x,y)$ are given
in a similar form.

(B) For \
$(\varphi^1_z)^2(x,y)=\zeta(x)\sin^2\varphi(x,y)-\eta(y)\cos^2\varphi(x,y)$
\ with {\rm analytic} functions $\zeta(x)$ and $\eta(y)$,
$\varphi(x,y)$ and $\varphi^1_z(x,y)$ are given, if there is
an {\rm analytic} function $\varphi(x,y)$ such that it
satisfies the following conditions (1) and (2): \

With
$${\hat A}:=-\frac{1}{2(\varphi^1_z)^2}(\zeta'\sin\varphi+2(\zeta+\eta)\varphi_x\cos\varphi) \hspace{0.4cm}{and}\hspace{0.4cm}
{\hat B}:=\frac{1}{2(\varphi^1_z)^2}(-\eta'\cos\varphi+2(\zeta+\eta)\varphi_y\sin\varphi),$$ 

(1) \ $(\zeta+\eta)\varphi_{xy}+\frac{1}{2}(\eta'\varphi_x+\zeta'\varphi_y)=-{\hat A}{\hat B}(\varphi^1_z)^2$ \ holds.

(2) \ There are functions $S=S(x,y)$, $T=T(x,y)$ such that \ $S_x=\varphi_x(L\varphi)$, \ $T_y=\varphi_y(L\varphi)$ \ and

 \ \ $L\varphi=S\cot\varphi-T\tan\varphi$.

Then, for any pair $\varphi(x,y)$ and $\varphi^1_z(x,y)$ in
the class (A) or (B), an analytic $2$-metric \ ${\hat
  g}:={\hat A}^2(dx)^2+{\hat B}^2(dy)^2$ \ with constant
Gauss curvature $-1$ is determined and a one-parameter
family $g^c$ of conformally flat 3-metrics given by (1) is obtained
via
evolution of orthogonal $2$-metrics with constant Gauss curvature $-1$ issuing from ${\hat g}$. 

Conversely, let ${\hat g}$ be an orthogonal analytic $2$-metric with
constant Gauss curvature $-1$. If there is a one-parameter family
$g^c$, $c\in R\setminus\{0\}$, of conformally flat 3-metrics given by (1) such that
their evolutions determined by $g^c$ satisfy \ ${\hat g^c}(0)={\hat g}$, \ then ${\hat g}$ is determined from some $\varphi(x,y)$ and $\varphi^1_z(x,y)$ in (A) or (B).
 }  \\

In this case, $g^c$ and $g^{c'}$ give distinct Guichard nets if \ $c\neq c'$ (Theorem 7 in \S4). 
 
The class (A) (resp. (B)) is characterised by the condition
on $\varphi(x,y)$ such that \
$(\varphi_{xy}-2\varphi_x\varphi_y\cot2\varphi)(x,y)=0$ \
(resp. $(\varphi_{xy}-2\varphi_x\varphi_y\cot2\varphi)(x,y)\neq
0$) \ (see Corollary 3.3 in \S3.2).
Main Theorem 2 proceeds by applying the Cauchy-Kovalevskaya
Theorem (which is why our data must be real-analytic) to
solve the system (2) with initial data at $z=0$.
For $\varphi(x,y)$ and $\varphi_z^1(x,y)$ in (A), respectively (B), 
we have \
$(L\psi^c)(x,y)=(1/2)[c^2\zeta(x)-\varphi_y^2/\cos^2\varphi]-\varphi_{yy}\tan\varphi$
\ and \
$(L\psi^c)(x,y)=(c^2/2)(\zeta(x)+\eta(y))+S(x,y)+T(x,y)$, \
respectively, and these equations determine the initial
$\psi^c(x,y)$ by solving a wave equation. From Main Theorem 2, we obtain many initial metrics ${\hat g}$ belonging to (A) by taking arbitrary $\zeta(x)$ and $D(y)$, and we shall also obtain many examples of ${\hat g}$ belonging to (B) (see \S2.2 and \S3.2). 

\vspace*{2mm}
Finally, remark that this analysis starts by distinguishing
the principal coordinate direction $z$.  However, a completely
analogous account may be given after distinguishing either the $x$- or
the $y$-direction
although, in these cases, the $2$-metrics will have indefinite
signature and constant curvature $1$.

\section{Existence condition for generic conformally flat
  hypersurfaces}
\label{sec:exist-cond-gener}

The existence of generic conformally flat hypersurfaces in 4-dimensional space forms is equivalent to that of functions $\varphi=\varphi(x,y,z)$ such that the following Riemannian 3-metric $g$ determined from $\varphi$ are conformally flat:
$$g=\cos^2\varphi dx^2+\sin^2\varphi dy^2+dz^2. \eqno{(1.1)}$$
Then, two conformally flat 3-metrics $g$ determined from $\varphi(x,y,z)$ and $\tilde{\varphi}(x,y,z)$ define the same Guichard net if and only if there are three constants $a_1$, $a_2$ and $a_3$ such that \ ${\tilde \varphi}(x,y,z)=\varphi(\pm x+a_1,\pm y+a_2,\pm z+a_3)$, \ as mentioned in the introduction.
That is, $\varphi$ is determined up to parameter shifts.
Furthermore,  
such a $3$-metric $g$ is conformally flat if and only if the covariant derivative $\nabla S$ of the Schouten tensor $S$ is totally symmetric, where \ $S=Ric-(R/4)g$ \ for the Ricci curvature $Ric$ and the scalar curvature $R$ of $g$. 
In terms of $\varphi$, the condition for $g$ to be conformally flat reads:\\

{\pro}1.1. {\it \ \ A metric $g$ given by (1.1) is conformally flat if and only if the function $\varphi$ satisfies the following four equations: 
$$ \varphi_{xyz}+\varphi_x\varphi_{yz}\tan\varphi-\varphi_y\varphi_{xz}\cot\varphi=0,  \leqno{(1)}$$
$$\frac{\varphi_{xxx}-\varphi_{yyx}+\varphi_{zzx}}{2}-\frac{(\varphi_{xx}-\varphi_{yy})\cos2\varphi-\varphi_{zz}}{\sin2\varphi}\varphi_x-\varphi_{xz}\varphi_z\cot\varphi=0,  \leqno{(2)}$$
$$\frac{\varphi_{xxy}-\varphi_{yyy}-\varphi_{zzy}}{2}-\frac{(\varphi_{xx}-\varphi_{yy})\cos2\varphi-\varphi_{zz}}{\sin2\varphi}\varphi_y-\varphi_{yz}\varphi_z\tan\varphi=0,  \leqno{(3)}$$
$$\frac{\varphi_{xxz}+\varphi_{yyz}+\varphi_{zzz}}{2}+\frac{\varphi_{xx}-\varphi_{yy}-\varphi_{zz}\cos2\varphi}{\sin2\varphi}\varphi_z-\varphi_x\varphi_{xz}\cot\varphi+\varphi_y\varphi_{yz}\tan\varphi=0.  \leqno{(4)}$$ } \\

\vspace*{1mm}
The four equations in Proposition 1.1 are equivalent to the fact that the following two differential 1-form $\alpha$ and 2-form $\beta$ determined from $\varphi$ are closed: 
$$\alpha=-\varphi_{xz}\cot\varphi dx+\varphi_{yz}\tan\varphi dy+\frac{\varphi_{xx}-\varphi_{yy}-\varphi_{zz}\cos2\varphi}{\sin2\varphi}dz,$$
$$\beta=\varphi_{xz}\cot\varphi dy\wedge dz +\varphi_{yz}\tan\varphi dz\wedge dx
-\frac{(\varphi_{xx}-\varphi_{yy})\cos2\varphi-\varphi_{zz}}{\sin2\varphi}dx\wedge dy.$$
More precisely, $\alpha$ is closed if and only if the first three equations (1)-(3) for $\varphi$ in Proposition 1.1 hold, and $\beta$ is closed if and only if the last equation (4) holds. 
Thus, the problem to find a generic conformally flat
hypersurface is reduced to that of finding a function $\varphi$ such
that the two differential forms $\alpha$ and $\beta$ are closed. 

From now on, let us assume that all functions are defined on a simply
connected domain $U=D\times I$ in $R^3=R^2\times R$, where $0\in I$.\\

{\thm}1. {\it \ \ For a given $\varphi(x,y,z)$ such that $d\alpha=d\beta=0$, there is a function $\psi(x,y,z)$ satisfying the following four equations: 
$$(1)\hspace{2mm} \psi_{xz}=-\varphi_{xz}\cot\varphi, \hspace{1.2cm} (2)\hspace{2mm}  \psi_{yz}=\varphi_{yz}\tan\varphi,$$
$$ (3)\hspace{2mm} \psi_{zz}=(\varphi_{xx}-\varphi_{yy})\sin2\varphi-(\psi_{xx}-\psi_{yy})\cos2\varphi,$$
$$ (4)\hspace{2mm} \varphi_{zz}=(\varphi_{xx}-\varphi_{yy})\cos2\varphi+(\psi_{xx}-\psi_{yy})\sin2\varphi.$$ 

Conversely, if there are two functions $\varphi$ and $\psi$
satisfying these four equations, then the 1-form $\alpha$
and 2-form $\beta$ determined by $\varphi$ are closed. 

In this case, we can assume that $\psi$ does not have any linear term for $x$, $y$, $z$.} \\

By Theorem 1, the system of the third order differential
equations for $\varphi$ in Proposition 1.1 are reduced to the
system of the second order differential equations for two
functions $\varphi$ and $\psi$.  However, $\psi$ is not
uniquely determined by $\varphi$ even if we insist on
vanishing linear term since, as we see in equations (3) and
(4), $\psi(x,y,z)$ has the ambiguity of terms $k(x+y)$ and
$\hat{k}(x-y)$ of 1-variable functions. We shall investigate
this fact in \S4, where we impose additional constraints (in
Proposition~4.1) after which $\psi$ is uniquely
determined by $\varphi$.

Theorem 1 is obtained from the following Proposition 1.2: \\

{\pro}1.2. {\it \ \ The existence of a function $\varphi(x,y,z)$ such that \ $d\alpha=d\beta=0$ \ is equivalent to the existence of  
functions $\varphi(x,y,z)$ and $\psi=\psi(x,y,z)$ such that $\varphi$ and $\psi$ satisfy the following four equations: 
$$(1)\hspace{2mm} \psi_{xz}=-\varphi_{xz}\cot\varphi, \hspace{2cm} (2)\hspace{2mm}  \psi_{yz}=\varphi_{yz}\tan\varphi,$$ 
$$ (3)\hspace{2mm}  \psi_{zz}=\frac{\varphi_{xx}-\varphi_{yy}-\varphi_{zz}\cos2\varphi}{\sin2\varphi}, \hspace{0.5cm} 
(4)\hspace{2mm} \psi_{xx}-\psi_{yy}=-\frac{(\varphi_{xx}-\varphi_{yy})\cos2\varphi-\varphi_{zz}}{\sin2\varphi}.$$ 
Then, we can choose the function $\psi$ such that it does not have any linear term for $x$, $y$, $z$. } \\

We can rewrite (3) and (4) in Proposition 1.2 to (3) and (4) in Theorem 1, in particular, Theorem 1-(3) is obtained by substituting $\varphi_{zz}$ in Proposition 1.2-(4) into (3).   \\

{\sc Proof of Proposition 1.2}. \ \ Firstly, we assume $d\alpha=0$.

$d\alpha=0$ is equivalent to the existence of a function $p=p(x,y,z)$ which satisfies the following three equations: 
$$p_x=-\varphi_{xz}\cot\varphi, \hspace{1cm} p_y=\varphi_{yz}\tan\varphi,\hspace{1cm} p_z=\frac{\varphi_{xx}-\varphi_{yy}-\varphi_{zz}\cos2\varphi}{\sin2\varphi}.\hspace{1.2cm}$$
Such a function $p$ is determined up to a constant term. 
We choose the constant as zero, then $p$ is uniquely determined from $\varphi$. We define a function $\hat{\psi}=\hat{\psi}(x,y,z)$ by 
$$\hat{\psi}(x,y,z):=\int^z_{0}p(x,y,z)dz.$$
Then $\hat{\psi}$ satisfies \ $\hat{\psi}(x,y,0)=0$ \ and
$$(1)\hspace{2mm} \hat{\psi}_{xz}=-\varphi_{xz}\cot\varphi \hspace{1cm} (2)\hspace{2mm}  \hat{\psi}_{yz}=\varphi_{yz}\tan\varphi \hspace{1cm}
(3)\hspace{2mm}  \hat{\psi}_{zz}=\frac{\varphi_{xx}-\varphi_{yy}-\varphi_{zz}\cos2\varphi}{\sin2\varphi}. $$
We note that, even if we replace $\hat{\psi}$ by \ $\psi(x,y,z)=\hat{\psi}(x,y,z)+f(x,y)$ \ with function $f(x,y)$, $\psi$ also satisfies the equations (1), (2), (3) and \ $\psi(x,y,0)=f(x,y)$, \ that is, \ $d\alpha=0$ \ determines $\psi_z$ (not $\psi$). 

Next, we express the 2-form $\beta$ by using $\hat{\psi}$ as follows: 
$$\beta=-\hat{\psi}_{xz}dy\wedge dz+\hat{\psi}_{yz}dz\wedge dx-\frac{(\varphi_{xx}-\varphi_{yy})\cos2\varphi-\varphi_{zz}}{\sin2\varphi}dx\wedge dy.$$
The condition $d\beta=0$ is equivalent to the equation
$$(\hat{\psi}_{xx}-\hat{\psi}_{yy})_z=-[\frac{(\varphi_{xx}-\varphi_{yy})\cos2\varphi-\varphi_{zz}}{\sin2\varphi}]_z.$$
Thus, there exists a function $\hat{f}(x,y)$ such that 
$$\hat{\psi}_{xx}-\hat{\psi}_{yy}+\hat{f}(x,y)=-\frac{(\varphi_{xx}-\varphi_{yy})\cos2\varphi-\varphi_{zz}}{\sin2\varphi}.$$
We find a function $f(x,y)$ by solving the wave equation $f_{xx}-f_{yy}=\hat{f}.$  
However, such a function $f(x,y)$ is not unique, i.e., we
can replace $f(x,y)$ by $f(x,y)+k(x,y)$ with any function $k(x,y)$ satisfying $k_{xx}-k_{yy}=0$. 
Here, we may assume that $f(x,y)$ does not have any linear term for $x$, $y$. Even under this assumption, $f(x,y)$ still has the ambiguity of terms $k(x+y)$ and $\hat{k}(x-y)$ of 1-variable functions. 

Since $\hat{\psi}$ vanishes on $z=0$, $f(x,y)$ satisfies
$$(f_{xx}-f_{yy})(x,y)=-\frac{(\varphi_{xx}-\varphi_{yy})\cos2\varphi-\varphi_{zz}}{\sin2\varphi}(x,y,0).\eqno{(1.2)}$$
We now define the required function $\psi=\psi(x,y,z)$ by
$$\psi(x,y,z):=\hat{\psi}(x,y,z)+f(x,y), \ \ \ (\psi(x,y,0)=f(x,y)).\eqno{(1.3)}$$
Then, we obtain (4) in the Proposition.

In particular, we can express $\alpha$ and  $\beta$ in terms of the function $\psi$ as follows: 
$$\alpha=d(\psi_z)=\psi_{xz}dx+\psi_{yz}dy+\psi_{zz}dz,$$
$$\beta=d(\psi_ydx+\psi_xdy)
=-\psi_{xz}dy\wedge dz+\psi_{yz}dz\wedge dx+(\psi_{xx}-\psi_{yy})dx\wedge dy.$$
This fact shows that the converse of the statement in the Proposition is also true. \hspace{\fill}$\Box$\\

For $d\beta=0$, we also have the following fact:\\

{\pro}1.3 (\cite{hs2}). {\it \ \ Suppose $d\alpha=0$. Then, 
$d\beta=0$ holds if and only if the following equation is satisfied: 
$$[\psi_{zz}]_z=[-\Delta\psi+\{(\varphi_x)^2+(\varphi_y)^2+(\varphi_z)^2\}]_z,$$
where $\Delta\psi=\left(\frac{\partial^2}{\partial x^2}+\frac{\partial^2}{\partial y^2}\right)\psi$.} \\

{\sc Proof}. \ \ We firstly note that $d\beta=0$ is equivalent to Proposition 1.1-(4). 
With respect to the coefficients of $\alpha$, Proposition 1.1-(4) is formulated as the following: 
$$(-\varphi_{xz}\cot\varphi)_x+(\varphi_{yz}\tan\varphi)_y+\left(\frac{\varphi_{xx}-\varphi_{yy}-\varphi_{zz}\cos 2\varphi}{\sin2\varphi}\right)_z=(\varphi_x^2+\varphi_y^2+\varphi_z^2)_z,$$
that is, \ 
$(\psi_{xz})_x+(\psi_{yz})_y+(\psi_{zz})_z=(\varphi_x^2+\varphi_y^2+\varphi_z^2)_z$ \ is satisfied under the condition $d\alpha=0$. 
The proposition now follows. \hspace{\fill}$\Box$\\

\section{Geometrical meaning of functions $\varphi$ and
  $\psi$.}
\label{sec:geom-mean-funct}

In this section, in particular, in \S2.1, we study a geometrical meaning of the equations (1) and (2) in Theorem 1 (resp. Proposition 1.2). 
In \S2.2. we study the converse proposition of the result in \S2.1. 

In \S2.1, we assume that $g$ given by (1.1) is conformally flat and that $\varphi$ satisfies \ $\varphi_{xz}\neq 0$ \ and \ $\varphi_{yz}\neq 0$. \ 
We recall that, in the case \ $\varphi_{xz}=\varphi_{yz}=0$, \ the metric $g$ determined by $\varphi$ leads to a generic conformally flat hypersurface either of product-type or with cyclic Guichard net.

\subsection{Evolution of metrics on surfaces with constant Gauss
  curvature $-1$}
\label{sec:flow-metr-surf}

Let us define the functions $\hat{A}(x,y,z)$ and $\hat{B}(x,y,z)$ from (1) and (2) in Theorem 1 by 
$$ \hat{A}:=-\frac{\varphi_{xz}}{\varphi_z\sin\varphi}=\frac{\psi_{xz}}{\varphi_z\cos\varphi}, \hspace{1cm} \hat{B}:=\frac{\varphi_{yz}}{\varphi_z\cos\varphi}=\frac{\psi_{yz}}{\varphi_z\sin\varphi}.$$
Then, we have the following Theorem: \\

{\thm}2. {\it \ \ Suppose that $\varphi(x,y,z)$ and
$\psi(x,y,z)$ satisfy the equations (1), (2) in Theorem 1.
Let $\hat{A}$ and $\hat{B}$ be defined as above. 
Then, for each $z$, the Riemannian $2$-metric $\hat{g}(z)$ on
the ($x,y$)-plane,
$$\hat{g}(z):=\hat{A}^2(x,y,z)dx^2+\hat{B}^2(x,y,z)dy^2,  \eqno{(2.1)}$$
has constant Gauss curvature $K_{\hat{g}(z)}\equiv -1$.} \\

{\sc Proof}. \ \ Firstly, we have the following equations from the definitions of $\hat{A}$ and $\hat{B}$: 
$$\varphi_{xz}=-\hat{A}\varphi_z\sin\varphi,\hspace{0.8cm} \varphi_{yz}=\hat{B}\varphi_z\cos\varphi,\hspace{0.8cm}
\psi_{xz}=\hat{A}\varphi_z\cos\varphi,\hspace{0.8cm} \psi_{yz}=\hat{B}\varphi_z\sin\varphi.$$
Then, by the integrability condition $(\varphi_{xz})_y=(\varphi_{yz})_x$, 
we have
$$(\hat{A}_y-\hat{B}\varphi_x)\sin\varphi+(\hat{B}_x+\hat{A}\varphi_y)\cos\varphi=0. \eqno{(2.2)}$$
By the integrability condition $(\psi_{xz})_y=(\psi_{yz})_x$, we have
$$(\hat{A}_y-\hat{B}\varphi_x)\cos\varphi-(\hat{B}_x+\hat{A}\varphi_y)\sin\varphi=-\hat{A}\hat{B}.\eqno{(2.3)}$$
When we substitute (2.2) into (2.3), we obtain 
$$ \hat{B}_x+\hat{A}\varphi_y=\hat{A}\hat{B}\sin\varphi,\hspace{2cm} \hat{A}_y-\hat{B}\varphi_x=-\hat{A}\hat{B}\cos\varphi. \eqno{(2.4)}$$

The integrability condition \ $(\varphi_x)_y=(\varphi_y)_x$ \ implies 
$$\left[\frac{\hat{B}_x}{\hat{A}}\right]_x+\left[\frac{\hat{A}_y}{\hat{B}}\right]_y=(\hat{B}\sin\varphi)_x-(\hat{A}\cos\varphi)_y\\
=(\hat{B}_x+\hat{A}\varphi_y)\sin\varphi-(\hat{A}_y-\hat{B}\varphi_x)\cos\varphi=\hat{A}\hat{B},$$
which shows $K_{\hat{g}}\equiv =-1$. 

In this construction, we note that, for each $z_0$, the
metric ${\hat g}(z_0)$ is defined so long as (1) and (2) of
Theorem 1, viewed as equations
on $\phi_{z}|_{z=z_{0}}$, $\psi_{z}|_{z=z_0}$ are satisfied along $z=z_0$. \hspace{\fill}$\Box$\\

In the proof of Theorem 2, we have obtained the following Corollary.\\

{\cor}2.1. {\it \ \ We have 
$$\varphi_x=\frac{\hat{A}_y}{\hat{B}}+\hat{A}\cos\varphi, \hspace{1cm}
\varphi_y=-\frac{\hat{B}_x}{\hat{A}}+\hat{B}\sin\varphi, \leqno{(a)}$$
$$(\log|\varphi_z|)_x=-\hat{A}\sin\varphi,\hspace{2cm} (\log|\varphi_z|)_y=\hat{B}\cos\varphi, \leqno{(b)}$$ 
$$\psi_{xz}=-\varphi_{xz}\cot\varphi, \hspace{2cm}\psi_{yz}=\varphi_{yz}\tan\varphi.\leqno{(c)}$$ }

Now, when we regard $\varphi_z(x,y, z)$ and $\psi_z(x,y,z)$ as 2-variable functions of $x$ and $y$ with parameter $z$, we also have the following Corollary of Theorem 2. \\
 
{\cor}2.2. \ \ {\it We have
$${\hat g}(z)=\frac{1}{\varphi_z^2(x,y,z)}\left\{(d\varphi_z)^2(x,y,z)+(d\psi_z)^2(x,y,z)\right\}.$$} \\

{\sc Proof}. \ \ We have the following two equations:
$$d\varphi_z=\varphi_{xz}dx+\varphi_{yz}dy=-{\hat A}\varphi_z\sin\varphi dx+{\hat B}\varphi_z\cos\varphi dy,$$
$$d\psi_z=\psi_{xz}dx+\psi_{yz}dy={\hat A}\varphi_z\cos\varphi dx+{\hat B}\varphi_z\sin\varphi dy.$$
Hence, we have \ 
$(d\varphi_z)^2+(d\psi_z)^2=\varphi_z^2({\hat A}^2dx^2+{\hat B}^2dy^2).$ \hspace{\fill}$\Box$\\

{\sc Remark}. \ \ When we define 
$${\bar A}:=-\frac{\varphi_{xz}}{\sin\varphi}=\frac{\psi_{xz}}{\cos\varphi} \hspace{0.5cm}{\rm and}\hspace{0.5cm} {\bar B}:=\frac{\varphi_{yz}}{\cos\varphi}=\frac{\psi_{yz}}{\sin\varphi},$$
a metric \ ${\bar g}(z):={\bar A}^2(x,y,z)(dx)^2+{\bar B}^2(x,y,z)(dy)^2$ \ is flat for each $z$. In this case, we also have a similar fact to Corollary 2.1 and, in particular, (b) is replaced by the following (b$'$): 
$$\varphi_{zx}=-\bar{A}\sin\varphi,\hspace{2cm} \varphi_{zy}=\bar{B}\cos\varphi. \leqno{(b')}$$
There is a crucial difference between (b) and (b$'$), and it is
essential for our study to consider metrics ${\hat g}(z)$ with
constant Gauss curvature $-1$ (see Theorem 3-(2) below and Theorem 5, Proposition 3.4
in \S3.2).

\subsection{Characterization of $2$-metrics with constant Gauss
  curvature $-1$}
\label{sec:chact-2-metr}

Let ${\hat g}(z)$ be an evolution of orthogonal $2$-metrics with constant Gauss curvature $-1$, given in Theorem 2. 
Then, for each \ $z=z_0$, \ ${\hat g}(z_0)$ has been defined from $\varphi(x,y,z_0)$, $\varphi_z(x,y,z_0)$ and $\psi_z(x,y,z_0)$. 
Here, we study the converse construction. 

Let ${\hat g}$ be a (local-)Riemannian $2$-metric of $C^{\infty}$ with constant Gauss curvature $-1$, defined by 
$${\hat g}:={\hat A}^2(x,y)(dx)^2+{\hat B}^2(x,y)(dy)^2. \eqno{(2.5)}$$ 
In the following Theorem 3, we show that three functions
$\varphi(x,y,0)$, $\varphi_z(x,y,0)$ and $\psi_z(x,y,0)$ are
determined from ${\hat g}$.  Our notation anticipates that,
in arguments to follow, $\varphi_z(x,y,0)$ and
$\psi_z(x,y,0)$ will be the $z$-derivatives on $z=0$ of
functions $\varphi(x,y,z)$ and $\psi(x,y,z)$. 
However, in Theorem 3, we do not assume the existence of
such extensions and work only with $\varphi(x,y,0)$, $\varphi_z(x,y,0)$ and $\psi_z(x,y,0)$. \\

{\thm}3. \ \ {\it Let a $2$-metric ${\hat g}$ given by (2.5)
  have constant Gauss curvature $-1$.
  Then:

\noindent
(1)  A function $\varphi(x,y,0)$ is well-defined by the
equations (a)
$$
\begin{array}{l}
\varphi_x(x,y,0):=\left({\hat{A}_y}/{\hat{B}}\right)(x,y)+\hat{A}(x,y)\cos\varphi(x,y,0), \\  
\varphi_y(x,y,0):=-\left({\hat{B}_x}/{\hat{A}}\right)(x,y)+\hat{B}(x,y)\sin\varphi(x,y,0),\\
\end{array}
\leqno{(a)}
$$
i.e., \ $(\varphi_x)_y(x,y,0)=(\varphi_y)_x(x,y,0)$ \ is satisfied. 
In particular, for any given $\lambda\in R$, \ $\varphi(x,y,0)$ satisfying \ $\varphi(0,0,0)=\lambda$ is uniquely determined. 

\noindent
(2)  Functions $\varphi_z(x,y,0)$ and $\psi_z(x,y,0)$ are also well-defined by the following equations (b) and (c), respectively: 
$$(\log|\varphi_z|)_x(x,y,0):=-\hat{A}(x,y)\sin\varphi(x,y,0),\hspace{0.5cm} (\log|\varphi_z|)_y(x,y,0):=\hat{B}(x,y)\cos\varphi(x,y,0), \leqno{ \ \ \ (b)}$$
$$\psi_{zx}(x,y,0):=-(\varphi_{zx}\cot\varphi)(x,y,0), \hspace{0.5cm} \psi_{zy}(x,y,0):=(\varphi_{zy}\tan\varphi)(x,y,0), \leqno{ \ \ \ (c)}$$  
i.e., $(\varphi_{zx})_y=(\varphi_{zy})_x$ \ and \ $(\psi_{zx})_y=(\psi_{zy})_x$ are satisfied. 
In particular, $\varphi_z(x,y,0)$ and $\psi_z(x,y,0)$ are determined up to the same constant multiple $c\neq 0$, if $\psi_z(x,y,0)$ has no constant term. 
Furthermore, \ $\psi_{zxy}(x,y,0)=(\varphi_{zx}\varphi_y+\varphi_x\varphi_{zy})(x,y,0)$ \ and \ $\varphi_{zxy}(x,y,0)=-(\varphi_x\psi_{zy}+\varphi_y\psi_{zx})(x,y,0)$ \ hold. }  \\

{\sc Proof}. \ \ The statement (1) is obtained by direct calculation from the assumption that ${\hat g}$ has constant Gauss curvature $-1$. 
Here, we only show the statement (2). 

By (a), we have 
$$({\hat A}_y-{\hat B}\varphi_x)\sin \varphi+({\hat B}_x+{\hat A}\varphi_y)\cos\varphi=0, \eqno{(2.6)}$$
$$({\hat A}_y-{\hat B}\varphi_x)\cos \varphi-({\hat B}_x+{\hat A}\varphi_y)\sin\varphi=-{\hat A}{\hat B}. \eqno{(2.7)}$$
We may define $\varphi_{z}(x,y,0)$ and $\psi_{z}(x,y,0)$ by
$$\varphi_{zx}:=-{\hat A}\varphi_z\sin\varphi, \hspace{0.5cm}\varphi_{zy}:={\hat B}\varphi_z\cos\varphi \hspace{0.8cm}{\rm and}\hspace{0.8cm}
\psi_{zx}:={\hat A}\varphi_z\cos\varphi, \hspace{0.5cm}\psi_{zy}:={\hat B}\varphi_z\sin\varphi,$$
respectively, as \ $(\varphi_{zx})_y=(\varphi_{zy})_x$ \ (resp. $(\psi_{zx})_y=(\psi_{zy})_x$) \ is satisfied by (2.6) (resp. (2.7)). 
Then, these definitions imply (b) and (c), respectively. 
Thus, we have shown that (b) and (c) are well-defined for $\varphi_z(x,y,0)$ and $\psi_z(x,y,0)$, respectively. 

The last two equations of (2) follow from (c) by \ $(\varphi_{zx})_y(x,y,0)=(\varphi_{zy})_x(x,y,0)$ \ and \ $(\psi_{zx})_y(x,y,0)=(\psi_{zy})_x(x,y,0)$, \ respectively. 

The theorem now follows. \hspace{\fill}$\Box$\\

\vspace{2mm}
In general, it seems difficult to solve the initial condition (a), (b), (c) from a metric ${\hat g}$ with constant Gauss curvature $-1$. Here, we study the problem for the hyperbolic $2$-metric on the upper half plane.\\

{\exa}1. {\it \ \ Let \ ${\hat g}=(dx^2+dy^2)/(y+b)^2$ \ with a constant $b(>0)$. Then, we obtain the following functions from ${\hat g}$: \ For the sake of simplicity, we denote $x+a$ ($a$: const.), $y+b$ and $\varphi(x,y,0)$ by $x$, $y$ and $\varphi$, respectively. 
$$\cos\varphi=\frac{x^2-y^2}{x^2+y^2}, \hspace{0.5cm} \sin\varphi=\frac{2xy}{x^2+y^2},\hspace{0.5cm} \left({\rm i.e.}, \ \varphi=\arctan(\frac{2xy}{x^2-y^2}),\hspace{0.5cm} \lambda=\arctan(\frac{2ab}{a^2-b^2})\right),$$
$$\varphi_x=-\frac{2y}{x^2+y^2}, \hspace{1cm} \varphi_y=\frac{2x}{x^2+y^2},\hspace{1cm} \varphi_z=\frac{cy}{x^2+y^2}, \hspace{1cm} \psi_z=\frac{-cx}{x^2+y^2} \ \ (c\neq 0: {\rm const.}).$$

For the study in \S3.2, we list other equations obtained in this case.  Let \ $\varphi_z(x,y,0)=\varphi^c_z(x,y,0):=c\varphi^1_z(x,y,0)$ \ and \ $\psi_z(x,y,0)=\psi^c_z(x,y,0):=c\psi^1_z(x,y,0)$ \ by Theorem 3-(2). We define functions \ $\zeta=\zeta(x):={1}/{4x^2}$, \ $S=S(x,y):={4y^2}/{(x^2+y^2)^2}$ \ and \ $T=T(x,y):=-{4x^2}/{(x^2+y^2)^2}+{1}/{x^2}.$ \ Then, we have  
$$S_x=\varphi_x(L\varphi), \ \ \ \ \ \ T_y=\varphi_y(L\varphi),$$  
$$(\varphi^1_z)^2=\zeta\sin^2\varphi, \hspace{1cm} L\varphi(:=\varphi_{xx}-\varphi_{yy})=\frac{8xy}{(x^2+y^2)^2}, \hspace{1cm} L\varphi=S\cot\varphi-T\tan\varphi.$$

Next, $\psi(x,y,0)$ is not determined from a metric ${\hat g}$ in Theorem 3. However, we can determine $\psi(x,y,0)$ for this metric under the assumption that the following equations (d) and (e) are satisfied: 
$$(d) \ \ \psi_{xy}=\varphi_x\varphi_y, \ \ \ \ \ 
 \ (e) \ \ -\Delta\psi+\varphi_x^2+\varphi_y^2+\varphi_z^2=(L\varphi)\sin2\varphi-(L\psi)\cos2\varphi.$$ 
\noindent
Furthermore, $\psi^c(x,y,0):=\psi(x,y,0)$ for each $c$ is uniquely determined 
$$\psi^c=\log(x^2+y^2)-(1+\frac{c^2}{8})\log x.$$ 
under the additional condition for $\psi^c(x,y,0)$ not to have linear
terms with respect to $x$ and $y$. The meaning of (d) and (e) becomes clear in the next section and these facts are verified in \S4 in a general situation.  
In particular, (d) and (e) in this case are given as follows: 
$$\psi^c_{xy}=\varphi_x\varphi_y =-\frac{4xy}{(x^2+y^2)^2}.$$ 
$$-\Delta\psi^c+\varphi_x^2+\varphi_y^2+(\varphi^c_z)^2=(L\varphi)\sin2\varphi-(L\psi^c)\cos2\varphi=-(1+\frac{c^2}{8})\frac{1}{x^2}+\frac{4}{x^2+y^2}+\frac{c^2y^2}{(x^2+y^2)^2}.$$ } \\

{\sc Proof}. \ \ Here, we only show that $\psi$ is uniquely
determined as above by (d), (e) and the additional
condition, as other functions are directly obtained from the definitions. 

We firstly have 
$$\psi=\log(x^2+y^2)+X(x)+Y(y)$$
with suitable functions $X$ and $Y$ of one variable, by (d).
Then, from (e), we obtain 
$$8X''x^2y^2+2Y''(x^2-y^2)^2=(8+c^2)y^2. \eqno{(2.8)}$$
Taking first and second derivatives of (2.8) with respect to $x$, we have 
$$\frac{(X''x^2)'}{x}y^2+Y''x^2=Y''y^2 \hspace{1cm}{\rm and}\hspace{1cm}
\frac{1}{x}\left(\frac{(X''x^2)'}{x}\right)'=-\frac{2Y''}{y^2}=c_1 \ ({\rm const.}). $$
Substituting \ $2Y''=-c_1y^2$ \ into (2.8), we have \ 
$8X''x^2-c_1(x^4-2x^2y^2+y^4)=8+c^2$. \ 
This equation implies \ $c_1=0$. \ Then, we have \ $X''=(1+c^2/8)/x^2$. Thus, $\psi$ has been determined for each $c$. \hspace{\fill}$\Box$\\

\section{Choice of initial data}
\label{sec:choice-initial-data}

We firstly study the integrability conditions on $\varphi_z$ and $\psi_z$ in Theorem 1, in \S3.1. 
Next, in \S3.2, we study the relation between the equations
(3), (4) of Theorem 1 and orthogonal $2$-metrics ${\hat g}$
with constant Gauss curvature $-1$. Through these studies,
we determine a class of initial data ${\hat g}$ for our
system of evolution equations (2) mentioned in the introduction.

\subsection{Integrability condition on $\varphi_z$ and
  $\psi_z$}
\label{sec:integr-cond-varph}

The following Theorem 4 and Proposition 3.1 are fundamental for our study. \\

{\thm}4. \ \ {\it Let $\varphi(x,y,z)$ and $\psi(x,y,z)$ satisfy all equations (1), (2), (3) and (4) in Theorem 1. 
Then, we have the following facts (1) and (2):

(1) The conditions of \ $(\psi_{zx})_y=(\psi_{zy})_x$ \ and \ $(\varphi_{zx})_y=(\varphi_{zy})_x$, \ respectively, are given by 
$$\varphi_{zxy}+\varphi_x\psi_{zy}+\varphi_y\psi_{zx}=0, \eqno{(3.1)}$$
$$\psi_{zxy}=\varphi_x\varphi_{zy}+\varphi_y\varphi_{zx}.  \eqno{(3.2)}$$ 

(2) The equations obtained from \
$(\psi_{xz})_z=(\psi_{zz})_x$ \ and \
$(\psi_{yz})_z=(\psi_{zz})_y$, \ respectively, are the same as
those obtained from \ $(\varphi_{xz})_z=(\varphi_{zz})_x$ \ and \ $(\varphi_{yz})_z=(\varphi_{zz})_y$. \ Furthermore, these equations imply that there are two 1-variable functions $k(x+y)$ and ${\hat k}(x-y)$ such that \ ${\tilde \psi}(x,y,z):=\psi(x,y,z)+k(x+y)+{\hat k}(x-y)$ \ satisfies the following (3.3) and (3.4):
$${\tilde \psi}_{xy}=\varphi_x\varphi_y,  \eqno{(3.3)}$$
$$(L\varphi)\sin 2\varphi-(L{\tilde \psi})\cos 2\varphi=
-\Delta{\tilde \psi}+(\varphi_x)^2+(\varphi_y)^2+(\varphi_z)^2.  \eqno{(3.4)}$$}\\

We note that ${\tilde \psi}(x,y,z)$ in Theorem 4-(2) also satisfies the all equations of Theorem 1, (3.1) and (3.2). 
Hence, the equation (3.4) means that the equation
$$ {\tilde \psi}_{zz}=(\varphi_{xx}-\varphi_{yy})\sin2\varphi-({\tilde \psi}_{xx}-{\tilde \psi}_{yy})\cos2\varphi=-\Delta{\tilde \psi}+(\varphi_x)^2+(\varphi_y)^2+(\varphi_z)^2 $$
is satisfied, by Theorem 1-(3). \\

{\sc Proof}. \ \ The statement (1) is obtained by direct calculation. In this proof, we only verify the statement (2), in particular, as the parameter $z$ varies on some interval, then, for the first statement of (2), we only study the equations induced from \ $(\psi_{xz})_z=(\psi_{zz})_x$ \ and  \ $(\psi_{yz})_z=(\psi_{zz})_y$, \ as we can obtain the equations from \ $(\varphi_{xz})_z=(\varphi_{zz})_x$ \ and \ $(\varphi_{yz})_z=(\varphi_{zz})_y$ \ in the same way. 

Before proceeding, we pause to consider that our goal in
\S3.2 is to view $\varphi$, $\varphi_{z}$, $\psi$, $\psi_z$
as initial data along a hypersurface $z=z_0$ and we want to
know under what conditions the conclusions of the present theorem hold
in that setting.  We shall therefore attempt to confine and
pinpoint our use of the equations of Theorem~1 and the
integrability conditions $(\psi_{zx})_z=(\psi_{zz})_x$ and
so on.

Let \ $L\varphi:=\varphi_{xx}-\varphi_{yy}$. \ Using $(\psi_{zx})_z=(\psi_{zz})_x$ and (1), (3), (4) of Theorem 1, we firstly have 
$$2\varphi_{xz}\varphi_z
=2\{(L\varphi)_x+2\varphi_x(L\psi)\}\sin\varphi\cos\varphi+2\{(L\psi)_x-2\varphi_x(L\varphi)\}\sin^2\varphi \eqno{(3.5)}$$
$$=[(L\varphi)\sin 2\varphi-(L\psi)\cos 2\varphi]_x+\{(L\psi)_x-2\varphi_x(L\varphi)\}. $$
For the second term of the last equation in (3.5), we have: 
$$(L\psi)_x-2\varphi_x(L\varphi)=[\Delta\psi-(\varphi_x^2+\varphi_y^2)]_x-2[\psi_{xy}-\varphi_x\varphi_y]_y.$$
By (3.5), we obtain the equation
$$[(L\varphi)\sin 2\varphi-(L\psi)\cos 2\varphi+\Delta{\psi}-(\varphi_x^2+\varphi_y^2+\varphi_z^2)]_x-2[{\psi}_{xy}-\varphi_x\varphi_y]_y=0. \eqno{(3.6)}$$

Similarly, using $(\psi_{yz})_z=(\psi_{zz})_y$ and (2), (3), (4) of Theorem 1, we have 
$$2\varphi_z\varphi_{yz}
=2[(L\varphi)_y+2\varphi_y(L\psi)]\sin\varphi\cos\varphi-2\{(L\psi)_y-2\varphi_y(L\varphi)\}\cos^2\varphi   \eqno{(3.7)}$$
$$=[(L\varphi)\sin2\varphi-(L\psi)\cos2\varphi]_y-\{(L\psi)_y-2\varphi_y(L\varphi)\}. $$
and
$$-\{(L\psi)_y-2\varphi_y(L\varphi)\}=[\Delta\psi-(\varphi_x^2+\varphi_y^2)]_y-2[\psi_{xy}-\varphi_x\varphi_y]_x.$$
Hence, we obtain
$$[(L\varphi)\sin 2\varphi-(L\psi)\cos 2\varphi+\Delta{\psi}-(\varphi_x^2+\varphi_y^2+\varphi_z^2)]_y-2[{\psi}_{xy}-\varphi_x\varphi_y]_x=0. \eqno{(3.8)}$$

Remark that the equivalence of (3.5), (3.7) with (3.6),
(3.8) uses only differentiations in $x,y$ and so is valid
along $z=z_0$.  

Furthermore, we have 
$$[(L\varphi)\sin 2\varphi-(L\psi)\cos 2\varphi+\Delta{\psi}-(\varphi_x^2+\varphi_y^2+\varphi_z^2)]_z=0$$
by Proposition 1.3. However, for this equation, our argument can not be restricted to \ $z=z_0$. 

The equations (3.6) and (3.8) imply that there are two 1-variable functions $l(x+y)$ and ${\hat l}(x-y)$ such that 
$$[(L\varphi)\sin 2\varphi-(L\psi)\cos 2\varphi+\Delta{\psi}-(\varphi_x^2+\varphi_y^2+\varphi_z^2)](x,y,z)=-l(x+y)-{\hat l}(x-y),$$
$$2[{\psi}_{xy}-\varphi_x\varphi_y](x,y,z)=-l(x+y)+{\hat l}(x-y),$$
as \ $p_x(x,y)=q_y(x,y)$ \ and \ $p_y(x,y)=q_x(x,y)$ \ imply \ $p_{xx}-p_{yy}=0$ \ and \ $q_{xx}-q_{yy}=0$. 

Finally, when we define \ $2k''(x+y):=l(x+y)$ \ and \ $2{\hat k}''(x-y):={\hat l}(x-y)$ \ and \ ${\tilde \psi}(x,y,z):={\psi}(x,y,z)+k(x+y)+{\hat k}(x-y)$, the function ${\tilde \psi}$ satisfies (3.3) and (3.4).

Remark again that this argument uses only (3.6) and (3.8)
and differentiations along $x,y$ and so hold on a fixed
coordinate surface $z=const$.

We have therefore proved the Theorem. \hspace{\fill}$\Box$\\

Equation (3.2) means that \ $\psi_{xy}(x,y,z)=(\varphi_x\varphi_y)(x,y,z)$ \ is satisfied for any $(x,y,z)$ if \ $\psi_{xy}(x,y,0)=(\varphi_x\varphi_y)(x,y,0)$ \ holds at any $(x,y,0)$.

In the following proposition, we give another proof of the fact that (3.5) and (3.7), respectively, are induced from \ $(\psi_{xz})_z=(\psi_{zz})_x$ \ and \ $(\psi_{yz})_z=(\psi_{zz})_y$, \ then it will be clear how (3.5), (3.7) are related with the equations in Proposition 1.1. Furthermore, we summarise equations equivalent to (3.5) and (3.7), which we have obtained in the proof of Theorem 4-(2). We shall use the result in \S3.2. \\

{\pro}3.1. \ \ {\it Suppose that all equations of Theorem 1 are satisfied. Then, the equation (3.9) below is satisfied for any $z$. 
Furthermore, suppose that all equations of Theorem 1 are satisfied at
arbitrarily fixed $z=z_0$.  Then, the following five statements (1), (2), (3), (4) and (5) are equivalent to each other at $z=z_0$. 

(1) \ The following equations are satisfied:
$$[{\psi}_{xz}+\varphi_{xz}\cot\varphi]_z=0, \hspace{1cm} [{\psi}_{yz}-\varphi_{yz}\tan\varphi]_z=0. \eqno{(3.9)}$$

(2) \ The following equations from (3.5) and (3.7) are satisfied:
$$\psi_{xzz}+(L\psi)_x-2\varphi_x(L\varphi)-2\varphi_z\varphi_{xz}=0,\hspace{0.8cm}
\psi_{yzz}-(L\psi)_y+2\varphi_y(L\varphi)-2\varphi_z\varphi_{yz}=0. $$

(3) \ The following equations from (3.5) and (3.7) are satisfied: 
$$(L\psi)_x=\frac{\varphi_{xz}\varphi_z}{\sin^2\varphi}-\{(L\varphi)_x+2\varphi_x(L\psi)\}\cot\varphi+2\varphi_x(L\varphi), $$
$$(L\psi)_y=-\frac{\varphi_{yz}\varphi_z}{\cos^2\varphi}+\{(L\varphi)_y+2\varphi_y(L\psi)\}\tan\varphi+2\varphi_y(L\varphi). $$ 

(4) \ (2) and (3) in Proposition 1.1 are satisfied:  
$$\varphi_{xzz}+(L\varphi)_x+2\varphi_x(L{\psi})+2\varphi_z{\psi}_{xz}=0,\hspace{0.8cm}
\varphi_{yzz}-(L\varphi)_y-2\varphi_y(L{\psi})+2\varphi_z{\psi}_{yz}=0.$$

(5) \ There are two 1-variable functions $k(x+y)$ and ${\hat k}(x-y)$ such that \ ${\tilde \psi}(x,y,z):=\psi(x,y,z)+k(x+y)+{\hat k}(x-y)$ \ satisfies (3.3) and (3.4): 
$${\tilde \psi}_{xy}=\varphi_x\varphi_y,  \hspace{0.8cm}
(L\varphi)\sin 2\varphi-(L{\tilde \psi})\cos 2\varphi=
-\Delta{\tilde \psi}+(\varphi_x)^2+(\varphi_y)^2+(\varphi_z)^2.  $$

In particular, the first (resp. second) equations of (1), (2), (3) and (4) are equivalent to each other.} \\

{\sc Proof}. \ \ It follows from (1) and (2) of Theorem 1 that (3.9) is satisfied for any $z$. 

From now on, let us fix $z=z_0$. Here, we only prove the equivalence between (1), (2) and (4) at $z=z_0$ simultaneously, as we showed other equivalences in the proof of Theorem 4.

Firstly, we study the equation 
$$0=({\psi}_{xz}\sin\varphi+\varphi_{xz}\cos\varphi)_z$$
$$=({\psi}_{xzz}-\varphi_z\varphi_{xz})\sin\varphi+(\varphi_{xzz}+\varphi_z{\psi}_{xz})\cos\varphi. \eqno{(3.10)}$$
When we substitute (3) and (4) of Theorem 1 into (3.10), we have 
$$0=({\psi}_{xz}\sin\varphi+\varphi_{xz}\cos\varphi)_z$$
$$=\{(L\varphi)_x+2\varphi_x(L{\psi})+\varphi_z{\psi}_{xz}\}\cos\varphi+\{(L{\psi})_x-2\varphi_x(L\varphi)-\varphi_z\varphi_{xz}\}\sin\varphi. \eqno{(3.11)}$$
From these equations, we have
$$0=(3.10)+(3.11)=$$
$$\{\varphi_{xzz}+(L\varphi)_x+2\varphi_x(L{\psi})+2\varphi_z{\psi}_{xz}\}\cos\varphi+\{{\psi}_{xzz}+(L{\psi})_x-2\varphi_x(L\varphi)-2\varphi_z\varphi_{xz}\}\sin\varphi. \eqno{(3.12)}$$
 
Now, in the equation (3.12), we have 
$${\rm the \ coefficient \ of} \ \sin\varphi=\tan\varphi\times({\rm the \ coefficient \ of} \ \cos\varphi), \eqno{(3.13)}$$
which shows that both sides of the equation (3.13) vanish. 

In fact, we consider the right hand side of (3.13):
$${\rm the \ coefficient \ of} \ \cos\varphi=2\left(\frac{(\varphi_{xx}-\varphi_{yy}+\varphi_{zz})_x}{2}+\varphi_xL({\psi})-\varphi_z\varphi_{xz}\cot\varphi\right), \eqno{(3.14)}$$
where $L({\psi})$ is given by Proposition 1.2-(4) from Theorem 1-(4), then the right hand side of (3.14) is same as the left hand side of Proposition 1.1-(2). 

Next, we shall prove the equality of (3.13): \ Substitute Theorem 1-(3) into ${\psi}_{zz}$ and express ${\psi}$ by $\varphi$, then we have 
$${\psi}_{xzz}+(L{\psi})_x-2\varphi_x(L\varphi)-2\varphi_z\varphi_{xz}$$
$$=(L{\psi})_x(1-\cos2\varphi)+(L\varphi)_x\sin2\varphi+2\varphi_x\varphi_{zz}-2\varphi_x(L\varphi)-2\varphi_z\varphi_{xz}$$
$$=\left[\frac{\varphi_{zz}-(L\varphi)\cos2\varphi}{\sin2\varphi}\right]_x(1-\cos2\varphi)+(L\varphi)_x\sin2\varphi+2\varphi_x(-(L\varphi)+\varphi_{zz})-2\varphi_z\varphi_{xz}$$
$$=\tan\varphi\left[(\varphi_{xx}-\varphi_{yy}+\varphi_{zz})_x-2\varphi_x\frac{(L\varphi)\cos2\varphi-\varphi_{zz}}{\sin2\varphi}-2\varphi_z\varphi_{xz}\cot\varphi\right],$$
which shows the equality of (3.13). 

Thus, we have that (3.10) holds if and only if  
$${\psi}_{xzz}+(L{\psi})_x-2\varphi_x(L\varphi)-2\varphi_z\varphi_{xz}=0, $$
which is (3.5) by the proof of Theorem 4. In consequence, the first equation of (1) is equivalent to the first equations of (2) and (4), respectively.

By starting from \ $[{\psi}_{yz}\cos\varphi-\varphi_{yz}\sin\varphi]_z(x,y,z)=0$, \ we also have that the equation is equivalent to  
$${\rm the \ left \ hand \ side \ of \ (3) \ in \ Proposition 1.1}={\psi}_{yzz}-(L{\psi})_y+2\varphi_y(L\varphi)-2\varphi_z\varphi_{yz}=0.$$
Hence, the second equation of (1) is equivalent to the second equations of (2) and (4), respectively. \hspace{\fill}$\Box$\\

We shall study more an interesting condition induced from (3.9), in the following section.

\subsection{Infinitesimal deformation of $2$-metrics with
  constant Gauss curvature $-1$}
\label{sec:infin-deform-2}

Let ${\hat g}$ be a (local-)Riemannian $2$-metric of $C^{\infty}$ with constant Gauss curvature $-1$, given by
$${\hat g}:={\hat A}^2(x,y)(dx)^2+{\hat B}^2(x,y)(dy)^2, $$
as in \S2.2. 
We now study the relation between such metrics ${\hat g}$ and the equations (3) and (4) of Theorem 1. 
In particular, we study an infinitesimal deformation of
${\hat g}$ in the $z$-direction of orthogonal metrics with constant Gauss curvature $-1$. 

For such a metric ${\hat g}$, we have obtained, in Theorem 3, functions $\varphi(x,y,0)$, $\varphi_z(x,y,0)$ and $\psi_z(x,y,0)$ satisfying 
$$\psi_{zx}(x,y,0)=-(\varphi_{zx}\cot\varphi)(x,y,0), \hspace{1cm} \psi_{zy}(x,y,0)=(\varphi_{zy}\tan\varphi)(x,y,0).   \eqno{(3.15)}$$
The system is uniquely determined by giving $\varphi(0,0,0)=\lambda$ and a constant $c\neq 0$, where we assumed that $\psi_z(x,y,0)$ has no constant term. 
We now formally assume the equations (3), (4) of Theorem 1 along $z=0$:   
$$
\begin{array}{l}
\psi_{zz}(x,y,0)=\left[(\varphi_{xx}-\varphi_{yy})\sin2\varphi-(L\psi)\cos2\varphi\right](x,y,0),\\
\varphi_{zz}(x,y,0)=\left[(\varphi_{xx}-\varphi_{yy})\cos2\varphi+(L\psi)\sin2\varphi\right](x,y,0)
\end{array}
\eqno{(3.16)}$$
with some function $(L\psi)(x,y,0)$ from which we will
recover  $\psi(x,y,0)$ by solving $(\psi_{xx}-\psi_{yy})=L\psi$. 

Under the preparation above, we recall the fact in Theorem 2 that the
existence of an evolution of orthogonal metrics $\hat{g}(z)$ with constant Gauss curvature $-1$ has been obtained from the equation
$$\psi_{zx}(x,y,z)=-(\varphi_{zx}\cot\varphi)(x,y,z), \hspace{1cm} \psi_{zy}(x,y,z)=(\varphi_{zy}\tan\varphi)(x,y,z). $$
Hence, for any fixed \ $z=z_0$, \ the condition for $\hat{g}(z_0)$ to deform infinitesimally in $z$-direction to orthogonal metrics with constant Gauss curvature $-1$, is given by the equations
$$\left[\psi_{zx}+\varphi_{zx}\cot\varphi \right]_z(x,y,z_0)=0, \hspace{1cm} \left[\psi_{zy}-\varphi_{zy}\tan\varphi\right]_z(x,y,z_0)=0.   $$

Applying the fact above, the condition for $\hat{g}$ to have infinitesimal deformation in $z$-direction to orthogonal metrics with constant Gauss curvature $-1$, is given by the equations
$$\left[\psi_{zx}+\varphi_{zx}\cot\varphi \right]_z(x,y,0)=0, \hspace{1cm} \left[\psi_{zy}-\varphi_{zy}\tan\varphi\right]_z(x,y,0)=0,   \eqno{(3.17)}$$
where we do not necessarily assume the existence of any
extensions of $\varphi(x,y,0)$, $\varphi_z(x,y,0)$ and
$\psi_z(x,y,0)$ around $z=0$, that is, we interpret
derivatives in (3.17) as \
$\psi_{zxz}(x,y,0):=\psi_{zzx}(x,y,0)$, \
$[\cot\varphi]_z(x,y,0):=-(\varphi_z/\sin^2\varphi)(x,y,0)$
and so on. Hence, (3.17) means that Proposition 3.1-(1) is
satisfied at $z=0$.  It then follows that statements
(2)--(5) of Proposition~3.1 hold along $z=0$.

Our aim here is to study the condition (3.17) for $\psi(x,y,0)$ only under the conditions (3.15) and (3.16) at \ $z=0$. \\

{\pro}3.2. \ \ {\it Let $\hat{g}$ be an orthogonal Riemannian $2$-metric with constant Gauss curvature $-1$, given as above. 
Let us take a system of functions $\varphi(x,y,0)$, $\varphi_z(x,y,0)$ and $\psi_z(x,y,0)$ determined from $\hat{g}$, by arbitrarily fixed $\lambda$ and $c$.  
Suppose that (3.16) and (3.17) are satisfied with some function \ $(L\psi)(x,y,0)$. \ Then, the following equation is satisfied:  
$$(L\psi)(x,y,0)\times(\varphi_{xy}\sin2\varphi-2\varphi_x\varphi_y\cos2\varphi)(x,y,0)=  \eqno{(3.18)}$$
$$\left[-\varphi_z\varphi_{zxy}+\varphi_{zx}\varphi_{zy}-((L\varphi)_{xy}+4\varphi_x\varphi_y(L\varphi))\frac{\sin2\varphi}{2}-\varphi_x(L\varphi)_y\sin^2\varphi+\varphi_y(L\varphi)_x\cos^2\varphi  \right](x,y,0).$$} \\

{\sc Proof}. \ \ We know that (3.17) is equivalent to
Proposition 3.1-(3) at $z=0$. We then arrive at (3.18) by direct calculation from \ $(L\psi)_{xy}(x,y,0)=(L\psi)_{yx}(x,y,0)$.  \hspace{\fill}$\Box$\\

Proposition~3.2 implies a necessary condition for $\hat{g}$ to
arise from a Guichard net: it is not necessarily the case
that $L\psi$ given by (3.18) actually satisfies Proposition
3.1-(3).  In general, this requirement amounts to a very
complicated differential equation for $\varphi(x,y,0)$ and
$\varphi_{z}(x,y,0)$.  However, we may simplify matters
somewhat by requiring solutions of (3.18) for all $c\neq 0$
as we now see.

We arbitrarily fix $\lambda$ such that \ $\varphi(0,0,0)=\lambda$ \ from now on: \ we wish to get conformally flat metrics with the Guichard condition (or conformally flat metrics given by (1.1)), then, for $\bar{\varphi}(x,y,z)$ such that \ $\bar{\varphi}(x,y,z):=\varphi(x+a,y+b,z)$ \ with constants $a$ and $b$, \ $\bar{\varphi}$ and $\varphi$ determine the same Guichard net. 
Hence, $\varphi(x,y,0)$ is uniquely determined from ${\hat
  g}$. However, $\varphi_z(x,y,0)$ depends on constants
$c\neq 0$ as well as ${\hat g}$ by Theorem 3. Let us denote $\varphi_z(x,y,0)=\varphi_z^c(x,y,0):=c\varphi_z^1(x,y,0)$. \ Then, we have the following Corollary of (3.18): \\

{\cor}3.3. \ \ {\it Let ${\hat g}$ be an orthogonal Riemannian $2$-metric with constant Gauss curvature $-1$. Let $\varphi(x,y,0)$ and $\varphi_z^c(x,y,0)$ for any $c\neq 0$ be functions determined from ${\hat g}$ as above. Then, $\varphi(x,y,0)$ satisfies one of the following two cases (A) and (B):  

(A) \ $(\varphi_{xy}\sin2\varphi-2\varphi_x\varphi_y\cos2\varphi)(x,y,0)=0$. \ Then, for each $c$ we have 
$$[-\varphi^c_z\varphi^c_{zxy}+\varphi^c_{zx}\varphi^c_{zy}-((L\varphi)_{xy}+4\varphi_x\varphi_y(L\varphi))\frac{\sin2\varphi}{2}-\varphi_x(L\varphi)_y\sin^2\varphi+\varphi_y(L\varphi)_x\cos^2\varphi  ](x,y,0)=0.$$

(B) \ $(\varphi_{xy}\sin2\varphi-2\varphi_x\varphi_y\cos2\varphi)(x,y,0)\neq 0$. \ Then, for each $c$, $(L\psi^c)(x,y,0)$ is uniquely determined by (3.18). }\\

{\sc Remark}. \ \ Case A has a pretty geometric
interpretation: the vanishing of
$(\varphi_{xy}\sin2\varphi-2\varphi_x\varphi_y\cos2\varphi)(x,y,0)$
is equivalent to the vanishing of
$(\ln\frac{\cos\varphi}{\sin\varphi})_{xy}(x,y,0)$ which happens
precisely when the coordinate surface $z=0$ is an isothermic
surface in any Guichard net $({R}^3,g)$ arising from
$\hat{g}$.
We thank the anonymous referee for this nice observation.\\

{\thm}5. \ \ {\it Let ${\hat g}$ be a $2$-metric with constant Gauss curvature $-1$. Suppose that $\varphi(x,y,0)$ and $\varphi^c_z(x,y,0):=c\varphi^1_z(x,y,0)$ determined by ${\hat g}$ satisfy the condition of Corollary 3.3-(A) for any $c\neq 0$ and that $\varphi(x,y,0)$, $\varphi^c_z(x,y,0)$ and $(L\psi^c)(x,y,0)$ satisfy Proposition 3.1-(3) at $z=0$ for any $c\neq 0$. Then, $\varphi(x,y,0)$ satisfies either \ $\cos^2\varphi(x,y,0)=1/(1+e^{D(y)}) \ or \ \cos^2\varphi(x,y,0)=1/(1+e^{C(x)})$, \ where $C(x)$ and $D(y)$ are any non-constant functions of one-variable. Furthermore, in the case of \ $\cos^2\varphi(x,y,0)=1/(1+e^{D(y)})$, \ we have 
$$(\varphi_z^c)^2=c^2\zeta(x)\sin^2\varphi, \hspace{1cm} L\psi^c=(1/2)[c^2\zeta(x)-\varphi_y^2/\cos^2\varphi]-\varphi_{yy}\tan\varphi,$$
where $\zeta(x)>0$ is any non-constant one-variable function. 

Conversely, if we define $\varphi(x,y,0)$, $\varphi_z^c(x,y,0)$ and $(L\psi^c)(x,y,0)$ for any $D(y)$ and $\zeta(x)>0$ as above, then an orthogonal $2$-metric ${\hat g}$ with constant Gauss curvature $-1$, which is independent of $c$, 
is determined such that $\varphi(x,y,0)$ and $\varphi^c_z(x,y,0)$ for ${\hat g}$ satisfy the condition of Corollary 3.3-(A) and that $\varphi(x,y,0)$, $\varphi^c_z(x,y,0)$ and $(L\psi^c)(x,y,0)$ satisfy Proposition 3.1-(3) at $z=0$.

In the case of \ $\cos^2\varphi(x,y,0)=1/(1+e^{C(x)})$, \ we also have similar results. } \\

We can assume that $\psi^c(x,y,0)$ determined from
$(L\psi^c)(x,y,0)$ in Theorem 5 satisfies Proposition
3.1-(5) at $z=0$, as the statements (1)-(5) at $z=0$ in
Proposition 3.1 are equivalent to each other. Hence, Theorem
5 provides many $2$-metrics ${\hat g}$ of this kind.\\
   
{\sc Proof}. \ \ Let ${\hat g}$ be a $2$-metric satisfying the assumption of the Theorem. 

We firstly consider the two equations in Corollary 3.3-(A). By the first equation, we have \ $\cos^2\varphi=1/(1+e^{(C(x)+D(y))})$ \ and \ $\sin^2\varphi=e^{(C(x)+D(y))}/(1+e^{(C(x)+D(y))})$, \ where $C(x)$ and $D(y)$ are one-variable functions. Since \ $[-\varphi^1_z\varphi^1_{zxy}+\varphi^1_{zx}\varphi^1_{zy}](x,y,0)=0$ \ by the first two terms in the left hand side of the second equation, we have \ $\varphi_z^1={\pm}e^{(F(x)+G(y))}$.  

Next, let \ $R(x,y,0,c^2):=(L\psi^c)(x,y,0)$ \ be a solution of Proposition 3.1-(3). 
Then, we have 
$$(\partial R/\partial c^2)_x=\varphi^1_{zx}\varphi^1_z/\sin^2\varphi-2\varphi_x(\partial R/\partial c^2)\cot\varphi, \hspace{0.5cm}
(\partial R/\partial c^2)_y=-\varphi^1_{zy}\varphi^1_z/\cos^2\varphi+2\varphi_y(\partial R/\partial c^2)\tan\varphi$$
by Proposition 3.1-(3). Hence, there are functions $\check{\zeta}(x,c^2)$, $\check{\eta}(y,c^2)$ such that 
$$(\partial R/\partial c^2)\sin^2\varphi=(1/2)[(\varphi^1_z)^2+\check{\eta}(y,c^2)], \hspace{0.5cm} (\partial R/\partial c^2)\cos^2\varphi=(1/2)[-(\varphi^1_z)^2+\check{\zeta}(x,c^2)],$$
and we have \ 
$(\varphi_z^1)^2=\check{\zeta}(x,c^2)\sin^2\varphi-\check{\eta}(y,c^2)\cos^2\varphi.$

Now, we have obtained
$$\frac{\check{\zeta}(x,c^2)e^{(C(x)+D(y))}-\check{\eta}(y,c^2)}{1+e^{(C(x)+D(y))}}=e^{2F(x)}e^{2G(y)}(=(\varphi_z^1)^2).$$
If \ $\chi(c^2):=\check{\zeta}(x,c^2)=-\check{\eta}(y,c^2)$, \ then \ $\chi(c^2)=e^{2(F(x)+G(y))}=(\varphi_z^1)^2.$ \ Since $\varphi^1_z$ is independent of $c^2$, $\chi(c^2)$ is constant and $F(x)$, $G(y)$ are also constants, which is contradiction to \ $\varphi^1_{zx}\neq 0$ \ and \ $\varphi^1_{zy}\neq 0$. \ Hence, this case does not occur. 

Otherwise, we use \
$(1+e^{(C(x)+D(y))})^{-1}=\Sigma_{n=0}^{\infty}(-e^{(C(x)+D(y))})^n$,
\ where we assumed \ $e^{(C(x)+D(y))}<1$ \ in the
neighborhood of $(0,0)$. If \ $e^{(C(x)+D(y))}>1$, \ then we
can replace \ $e^{(C(x)+D(y))}<1$ \ by \
$\{e^{(C(x)+D(y))}-a \}/(1+a)<1$ \ with a suitable constant
$a$ from \
$1+e^{(C(x)+D(y))}=(1+a)[1+\{e^{(C(x)+D(y))}-a\}/(1+a)]$. \
Then, we have at least \ $C(x)=0$ \ or \ $D(y)=0$, and may
assume \ $C(x)=0$. \ Indeed, in the case of \ $D(y)=0$, \
the argument below proceeds in the same way when we consider
\
$(\varphi_z^1)^2=[\check{\zeta}(x,c^2)-\check{\eta}(y,c^2)e^{-C(x)}]/(1+e^{-C(x)})=e^{2F(x)}e^{2G(y)}$.

Now, let us assume \ $C(x)=0$. \ Since \ ${\check \zeta}(x,c^2)e^{D(y)}-{\check \eta}(y,c^2)=[{\check \zeta}(x,c^2)-{\check \eta}(y,c^2)/e^{D(y)}]e^{D(y)}$, \ we have \ ${\check \eta}(y,c^2)=h(c^2)e^{D(y)}$ \ and that $e^{D(y)}$ really depends on $y$ since \ $G'(y)\neq 0$. \ 
We also obtain \ ${\check \zeta}(x,c^2)-h(c^2)=\zeta(x)$, \ where $\zeta(x)$ is independent of $c^2$ from \ $(\varphi_z^1)^2=e^{2(F(x)+G(y))}$. \ In consequence, we have 
$$(\varphi_z^1)^2=\zeta(x)\sin^2\varphi, \hspace{1cm} \partial R/\partial c^2=(1/2)[\zeta(x)+h(c^2)/\cos^2\varphi],$$
that is,
$$(\varphi^c_z)^2=c^2\zeta(x)\sin^2\varphi, \hspace{1cm} L\psi^c=(1/2)[c^2\zeta(x)+H(c^2)/\cos^2\varphi]+I(x,y),$$
where $H'(c^2)=h(c^2)$ and that $I(x,y)$ is independent of $c^2$. 

On the other hand, we consider the equations of Proposition 3.1-(3) under the condition \ $C(x)=0$, \ i.e., \ $\varphi_x(x,y,0)=0$ \ and \ $L\varphi=-\varphi_{yy}$. \ Then, there are functions ${\tilde \zeta}(x,c^2)$ and ${\tilde \eta}(y,c^2)$ such that \ $(L\psi^c)\sin^2\varphi=(1/2)((\varphi_z^c)^2+{\tilde \eta}(y,c^2))$ \ and \ 
$(L\psi^c)\cos^2\varphi=-(1/2)((\varphi_z^c)^2+\varphi_y^2-{\tilde \zeta}(x,c^2))-\varphi_{yy}\sin\varphi\cos\varphi$. \ Hence, we have 
$$L\psi^c=-(1/2)\varphi_y^2-\varphi_{yy}\sin\varphi\cos\varphi+(1/2)({\tilde \zeta}(x,c^2)+{\tilde \eta}(y,c^2)),$$
$$(\varphi_z^c)^2={\tilde \zeta}(x,c^2)\sin^2\varphi-{\tilde \eta}(y,c^2)\cos^2\varphi-\varphi_y^2\sin^2\varphi-2\varphi_{yy}\sin^3\varphi\cos\varphi.$$
Then, we have \ ${\tilde \zeta}(x,c^2)=c^2\zeta(x)$, \ ${\tilde \eta}(y,c^2)=:\eta(y)$ \ and \ $\eta(y)\cos^2\varphi=-\varphi_y^2\sin^2\varphi-2\varphi_{yy}\sin^3\varphi\cos\varphi$ \ by $(\varphi_z^c)^2$. Furthermore, we have \ $h(c^2)=H(c^2)=0$ \ by $L\psi^c$ and \ ${\tilde \eta}(y,c^2)=\eta(y)$. 

By the argument above, we obtain, with \ $\cos^2\varphi=1/(1+e^{D(y)})$,
$$(\varphi_z^c)^2=c^2\zeta(x)\sin^2\varphi,\hspace{0.5cm} L\psi^c=-(1/2)\varphi_y^2-\varphi_{yy}\sin\varphi\cos\varphi+(1/2)(c^2\zeta(x)+\eta(y)),$$
where \ $\eta(y):=[-\varphi_y^2\sin^2\varphi-2\varphi_{yy}\sin^3\varphi\cos\varphi]/\cos^2\varphi$ \ and that $D(y)$ and $\zeta(x)>0$ can be taken arbitrarily. These are functions $\varphi^c_z(x,y,0)$ and $(L\psi^c)(x,y,0)$ in the Theorem. 

Conversely, these functions satisfy Proposition 3.1-(3) and determine $2$-metrics ${\hat g}$ with constant Gauss curvature $-1$ by Theorem 2, as there is a function $\psi^c_z(x,y,0)$ such that \ $\psi^c_{zx}=-\varphi^c_{zx}\cot\varphi$ \ and \ $\psi^c_{zy}=\varphi^c_{zy}\tan\varphi$ \ for each pair of $\varphi$ and $\varphi^c_z$. Furthermore, these functions $\varphi(x,y,0)$, $\varphi^c_z(x,y,0)$ and $\psi^c_z(x,y,0)$ are also defined from such a ${\hat g}$, by Theorem 3. 

We can also obtain similar results in the case \ $\cos^2\varphi=1/(1+e^{C(x)})$. \hspace{\fill}$\Box$\\

Next, we study the condition on ${\hat g}$ in the case of Corollary 3.3-(B) such that \ $(L\psi^c)_{xy}(x,y,0)=(L\psi^c)_{yx}(x,y,0)$. Then, $(L\psi^c)(x,y,0)$ is divided into two terms by the expression (3.18): \ 
$$(L\psi^c)(x,y,0)=c^2P(x,y)+Q(x,y),$$  
where 
$$P(x,y):=\left(\frac{d}{dc^2}L\psi^c\right)(x,y,0) \ \left(=\frac{-\varphi^1_z\varphi^1_{zxy}+\varphi^1_{zx}\varphi^1_{zy}}{\varphi_{xy}\sin2\varphi-2\varphi_x\varphi_y\cos2\varphi}(x,y,0)\right),$$ 
$$Q(x,y):=\frac{-((L\varphi)_{xy}+4\varphi_x\varphi_y(L\varphi))\frac{\sin2\varphi}{2}-\varphi_x(L\varphi)_y\sin^2\varphi+\varphi_y(L\varphi)_x\cos^2\varphi}{\varphi_{xy}\sin2\varphi-2\varphi_x\varphi_y\cos2\varphi}(x,y,0).$$
Our assumption for $(L\psi^c)(x,y,0)$ of a $2$-metric ${\hat g}$ that Proposition 3.1-(3) is satisfied for arbitrary $c\neq 0$ is equivalent to the following equations at $z=0$:  
$$P_x=\frac{\varphi^1_{xz}\varphi^1_z}{\sin^2\varphi}-2\varphi_xP\cot\varphi, \hspace{1cm}
P_y=-\frac{\varphi^1_{yz}\varphi^1_z}{\cos^2\varphi}+2\varphi_yP\tan\varphi, \eqno{(3.19)}$$
$$Q_x=-\{(L\varphi)_x+2\varphi_xQ\}\cot\varphi+2\varphi_x(L\varphi),\hspace{0.5cm}
Q_y=\{(L\varphi)_y+2\varphi_yQ\}\tan\varphi+2\varphi_y(L\varphi). \eqno{(3.20)}$$\\

{\pro}3.4. \ \ {\it Let ${\hat g}$ be an orthogonal
$2$-metric with constant Gauss curvature $-1$. Let us define
$\varphi(x,y,0)$, $\varphi^c_z(x,y,0)$, $(L\psi^c)(x,y,0)$, $P(x,y)$
and $Q(x,y)$ for ${\hat g}$ as above under the assumption that
$(L\psi^c)(x,y,0)$ is expressed by (3.18). Suppose that $(L\psi^c)(x,y,0)$ with arbitrary $c\neq 0$ satisfies Proposition 3.1-(3) at $z=0$. Then, we have the following facts: 

\ (1) \ There are functions $\zeta=\zeta(x)$, $\eta=\eta(y)$ such that \ $(\varphi_z^1)^2=\zeta\sin^2\varphi-\eta\cos^2\varphi$ \ and \ $P=(\zeta+\eta)/2$.

\ (2) \ There are functions $S=S(x,y)$ and $T=T(x,y)$ such that \ $S_x=\varphi_x(L\varphi)$, \ $T_y=\varphi_y(L\varphi)$, \ $L\varphi(:=\varphi_{xx}-\varphi_{yy})=S\cot\varphi-T\tan\varphi$ \ and \ $Q=S+T$.

\ (3) \ $(L\psi^c)(x,y,0)=c^2(\zeta(x)+\eta(y))/2+S(x,y)+T(x,y)$ \ is satisfied. 

Conversely, suppose that, for $\varphi(x,y,0)$ and $\varphi^1_z(x,y,0)$ determined from ${\hat g}$, there are functions $\zeta(x)$, $\eta(y)$, $S(x,y)$ and $T(x,y)$ satisfying (1) and (2). Then, if we take $(L\psi^c)(x,y,0)$ given in (3), $(L\psi^c)(x,y,0)$ satisfies Proposition 3.1-(3) at $z=0$, that is, for such a $2$-metric ${\hat g}$, $(L\psi^c)(x,y,0)$ is determined such that it satisfies Proposition 3.1-(3) and Corollary 3.3-(B).}\\

{\sc Proof}. \ Let us assume that $(L\psi^c)(x,y,0)$ is
given by (3.18) and $(L\psi^c)(x,y,0)$ with any $c\neq 0$ satisfies the equations of Proposition 3.1-(3), that is, $P(x,y)$ and $Q(x,y)$ satisfy (3.19) and (3.20), respectively. We shall verify that the assumption is equivalent to (1) and (2). 

Now, since we have the following equations from (3.19): 
$$[P\sin^2\varphi-(\varphi_z^1)^2/2]_x=[P\cos^2\varphi+(\varphi_z^1)^2/2]_y=0,$$there are functions $\zeta=\zeta(x)$ and $\eta=\eta(y)$ such that  
$$P\sin^2\varphi-(\varphi_z^1)^2/2=\eta/2, \hspace{1cm} P\cos^2\varphi+(\varphi_z^1)^2/2=\zeta/2.$$
Hence, we obtain \ $P=(\zeta+\eta)/2$ \ and $(\varphi^1_z)^2=\zeta\sin^2\varphi-\eta\cos^2\varphi$. 

Next, since we have the following equations from (3.20): 
$$[Q\sin^2\varphi+(L\varphi)\sin\varphi\cos\varphi]_x=\varphi_x(L\varphi), \hspace{1cm} [Q\cos^2\varphi-(L\varphi)\sin\varphi\cos\varphi]_y=\varphi_y(L\varphi),$$
there are functions $S=S(x,y)$ and $T=T(x,y)$ such that \ $S_x=\varphi_x(L\varphi)$, \ $T_y=\varphi_y(L\varphi)$, \  
$$Q\sin^2\varphi+(L\varphi)\sin\varphi\cos\varphi=S \hspace{0.5cm}{\rm and}\hspace{0.5cm} Q\cos^2\varphi-(L\varphi)\sin\varphi\cos\varphi=T$$
are satisfied. 
Hence, we obtain \ $Q=S+T$ \ and \ $L\varphi=S\cot\varphi-T\tan\varphi.$ 

In each argument above, the converse is also valid. Finally, we obtain $L\psi^c$ from \ $L\psi^c=c^2P+Q$. 

We note about the converse statement: (3.18) has been obtained from the assumption that $(L\psi^c)(x,y,0)$ satisfies Proposition 3.1-(3) (resp. Proposition 3.1-(1)). 
Furthermore, suppose that there is a solution $\varphi(x,y,0)$ and $\varphi_z(x,y,0)$ such that \ $\varphi_x(x,y,0)=\eta(y)=0$ \ in this case. \ Then, $L\psi^c=(1/2)[c^2\zeta(x)-(\varphi_y)^2]+const.$ \ is different from the ones in Theorem 5.  
This fact implies that there is not such a solution in this case. 
Hence, the $2$-metrics ${\hat g}$ obtained here are included in Corollary 3.3-(B).  \hspace{\fill}$\Box$\\

Now, $\varphi_z^c(x,y,0)$ has been determined from ${\hat g}$, by
Theorem 3. Hence, the property of $\varphi^1_z(x,y,0)$ in Proposition
3.4-(1) induces a condition for ${\hat g}$. Next, we study this condition. 

Let us assume \ $(\varphi_z^1)^2=\zeta\sin^2\varphi-\eta\cos^2\varphi$ \ as in Proposition 3.4. \ Then, for ${\hat g}={\hat A}^2(dx)^2+{\hat B}^2(dy)^2$, we have 
$${\hat A}=-\frac{1}{2(\varphi_z^1)^2}(\zeta'\sin\varphi+2(\zeta+\eta)\varphi_x\cos\varphi), \hspace{0.5cm} {\hat B}=\frac{1}{2(\varphi_z^1)^2}(-\eta'\cos\varphi+2(\zeta+\eta)\varphi_y\sin\varphi) \eqno{(3.21)}$$
by Theorem 2 and Corollary 2.1-(b). Furthermore, the condition that ${\hat g}$ has constant Gauss curvature $-1$ is equivalent to the existence of $\psi^1_z$ such that $\psi^1_{zx}=-\varphi^1_{zx}\cot\varphi$ and $\psi^1_{zy}=\varphi^1_{zy}\tan\varphi$, by Theorem 2, Corollary 2.1 and Theorem 3. By the integrability condition of $\psi^1_z$, we have the following Proposition.\\

{\pro}3.5. \ \ {\it A $2$-metric ${\hat g}={\hat A}^2(dx)^2+{\hat B}^2(dy)^2$ defined by (3.21) from \ $(\varphi_z^1)^2=\zeta\sin^2\varphi-\eta\cos^2\varphi$ \ with $\zeta(x)$ and $\eta(y)$ has the constant Gauss curvature $-1$, if and only if the following equation is satisfied: 
$$(\zeta+\eta)\varphi_{xy}+\frac{1}{2}(\eta'\varphi_x+\zeta'\varphi_y)=-{\hat A}{\hat B}(\varphi_z^1)^2.$$}\\

We have the following Theorem by summarising Propositions 3.4 and 3.5:  \\

{\thm}6. \ \ {\it For functions $\zeta=\zeta(x)$, $\eta=\eta(y)$ of one variable, let us set \ $(\varphi_z^1)^2(x,y,0):=(\zeta\sin^2\varphi-\eta\cos^2\varphi)(x,y,0)$. \ Suppose that there is a function $\varphi(x,y,0)$ such that it satisfies the following equations (1) and (2): 

\ (1) \ $(\zeta+\eta)\varphi_{xy}+\frac{1}{2}(\eta'\varphi_x+\zeta'\varphi_y)=-{\hat A}{\hat B}(\varphi_z^1)^2$, where ${\hat A}$ and ${\hat B}$ are given by (3.21). 

\ (2) \ There are functions $S=S(x,y)$, $T=T(x,y)$ satisfying \ $S_x=\varphi_x(L\varphi)$, \ $T_y=\varphi_y(L\varphi)$ \ and

$L\varphi=S\cot\varphi-T\tan\varphi.$

\noindent
Then, a $2$-metric \ ${\hat g}:={\hat A}^2(dx)^2+{\hat B}^2(dy)^2$ \ with constant Gauss curvature $-1$ and functions \ $(L\psi^c)(x,y,0):=(c^2/2)(\zeta +\eta)+S+T$, \ $\varphi^c_z(x,y,0):=c\varphi^1_z(x,y,0)$, \ $\psi^c(x,y,0)$ and $\psi^c_z(x,y,0)$ are determined.

Furthermore, let us  define $\varphi^c_{zz}$, $\psi^c_{zz}$ by (3.16). Then, we can choose a suitable $\psi^c(x,y,0)$ such that the system $\{\varphi,\psi^c,\varphi^c_z, \psi^c_{zz}\}$ of functions with arbitrary $c\neq 0$ satisfies Proposition 3.1-(5) at $z=0$. 

Conversely, if every one-parameter system $\{\varphi,\psi^c,\varphi^c_z, \psi^c_{zz}\}$ at $z=0$ for any $c\neq 0$ determined by a metric ${\hat g}$ with constant Gauss curvature $-1$ satisfies Proposition 3.1-(5) and \ $(\varphi_{xy}\sin2\varphi-2\varphi_x\varphi_y\cos2\varphi)(x,y,0)\neq 0$, \ then the metric ${\hat g}$ is obtained from $\varphi(x,y,0)$ and $\varphi^1_z(x,y,0)$ satisfying (1) and (2).}\\

{\sc Proof}. \ \ The condition (1) determines a $2$-metric
${\hat g}$ with constant Gauss curvature $-1$, as in
Proposition 3.5. $\varphi(x,y,0)$, $\varphi^c_z(x,y,0)$ and
$\psi^c_z(x,y,0)$ arise from the metric ${\hat g}$, by
Theorem 3. Then, these functions coincide with the ones stated in the Theorem by the construction of ${\hat g}$ in (3.21), Theorem 2, Corollary 2.1 and Theorem 3. 

Let \ $P:=(\zeta+\eta)/2$, \ $Q:=S+T$ \ and \ $L\psi^c=c^2P+Q$. \ For the $L\psi^c$, we define $\varphi^c_{zz}$ and $\psi^c_{zz}$ by (3.16). Then, $L\psi^c$ satisfies the equations of Proposition 3.1-(3) at $z=0$, by Proposition 3.4. 

Furthermore, since $L\psi^c=\psi^c_{xx}-\psi^c_{yy}$, we can  determine $\psi^c(x,y,0)$ up to two 1-variable functions $k(x+y)$ and ${\hat k}(x-y)$. Taking a suitable $\psi^c(x,y,0)$, the system $\{\varphi,\psi^c,\varphi^c_z, \psi^c_{zz}\}$ of functions satisfies Proposition 3.1-(5) at $z=0$. 

The converse also follows from Propositions 3.2, 3.4 and 3.5, as \ $(\varphi_{xy}\sin2\varphi-2\varphi_x\varphi_y\cos2\varphi)(x,y,0)\neq 0$ \ is the condition that ${\hat g}$ belongs to the case of Corollary 3.3-(B). \hspace{\fill}$\Box$\\

We study some examples of $\varphi(x,y,0)$ and $\varphi^1_z(x,y,0)$ in Theorem 6 (see Example 1 in \S.2.2 and Examples 3, 4 below). 

\vspace{2mm}
Now, let $M$ be the space of (local) orthogonal $2$-metrics ${\hat g}$ on $(x,y)$-plane with constant Gauss curvature $-1$. 
Let ${\hat g}$ be a metric of $M$ given in Theorem 5 or obtained by the procedure in Theorem 6. 
Then, ${\hat g}$ has a $z$-direction such that, if there is a curve
through ${\hat g}$ in $M$ which determines a conformally flat metric
$g$ with the Guichard condition, then the curve evolves in the direction at ${\hat g}$. Its direction is actually determined by a pair of $\varphi^c_{zz}(x,y,0)$ and $\psi^c_{zz}(x,y,0)$ (see Theorem 7 in \S4). 
In particular, the $z$-direction at ${\hat g}$ is determined by a 1-parameter family with parameter $c\neq 0$. 

We shall show in \S4 that such an {\it analytic} metric ${\hat g}$
really extends to an evolution of $2$-metrics ${\hat g}(z)$ for each $c\neq 0$, which determines a conformally flat metric $g^c$ with the Guichard condition. Then, $g^c$ and $g^{c'}$ have different conformal structures if $c\neq c'$ by the definition. To find generic conformally flat hypersurfaces was the problem to obtain general solutions $\varphi(x,y,z)$ of four complicated differential equations of third order in Proposition 1.1. 
In consequence, under a generic condition, the problem is
reduced to find functions $\varphi(x,y,0)$ and
$\varphi^1_z(x,y,0)$ stated in Theorem 6, as their functions
in Theorem 5 are already obtained explicitly. Here, we used
the term "generic'' in the meaning that ${\hat g}$ gives
rise to a one parameter family $g^c$. 

We note that the conditions (c) \ $\psi_{xy}(x,y,0)=(\varphi_x\varphi_y)(x,y,0)$ \ and 
$$\psi_{zz}(x,y,0)=[(L\varphi)\sin2\varphi-(L\psi)\cos2\varphi](x,y,0) \leqno{(d)}$$
$$=[-\Delta\psi+(\varphi_x)^2+(\varphi_y)^2+(\varphi_z)^2](x,y,0)$$
are satisfied, for metrics ${\hat g}$ given in Theorem 5 and obtained by the procedure in Theorem 6, i.e., Proposition 3.1-(5) is satisfied for such metrics ${\hat g}$. \\

{\exa}2. (Counter example) \ \ On $z=0$, we set 
$$\varphi(x,y):=x+y, \ \ \ \varphi^c_z(x,y):=\frac{ce^{y-x}}{2}(\cos\varphi+\sin\varphi), \ \ \ \psi_z^c:=\frac{ce^{y-x}}{2}(-\cos\varphi+\sin\varphi)$$
with constant $c\neq 0$. Then, for (b) and (c) of Theorem 3, we have
$$\psi^c_{zx}(x,y)=\varphi_{zy}^c(x,y)=ce^{y-x}\cos\varphi(x,y), \ \ \ \ \psi^c_{zy}(x,y)=-\varphi_{zx}^c(x,y)=ce^{y-x}\sin\varphi(x,y).$$
The $2$-metric \ ${\hat g}=4/(\cos\varphi+\sin\varphi)^2 \ ((dx)^2+(dy)^2)$ \ defined by the functions above has the constant Gauss curvature $-1$. Then, we obtain 
$$(L\psi^c)(x,y)=-\frac{c^2e^{2(y-x)}}{4\cos2\varphi(x,y)}$$ 
from (3.18) for ${\hat g}$. However, this $(L\psi^c)(x,y)$ does not
satisfy Proposition 3.1-(3). Hence, this metric ${\hat g}$ does not
extend into the $z$-direction.

In fact, from the first equation of Proposition 3.1-(3) we have 
$$\sin\varphi+\cos2\varphi(\cos\varphi+\sin\varphi)=\frac{\cos\varphi}{\cos2\varphi},$$
and from the second equation, we have
$$\cos\varphi-\cos2\varphi(\cos\varphi+\sin\varphi)=-\frac{\sin\varphi}{\cos2\varphi}.$$
If these two equations are satisfied, then we obtain 
$$\cos2\varphi(\cos\varphi+\sin\varphi)=\cos\varphi-\sin\varphi$$
by adding two equations. 
Then, we simultaneously have \ $\cos2\varphi=\pm 1$, 
which can not occur. \\

{\exa}3. \ \ Let us take \
$(\varphi^1_z)^2=c_1\sin^2\varphi-c_2\cos^2\varphi$, \ that
is, \ $\zeta(x)=c_1$ \ and \ $\eta(y)=c_2$. \ Then, the
function $\varphi(x,y,0)$ such that \
$\varphi_x=c_3\sqrt{c_1\sin^2\varphi-c_2\cos^2\varphi}$, \
$\varphi_y=c_4\sqrt{c_1\sin^2\varphi-c_2\cos^2\varphi}$
satisfies the condition (1) and (2) in Theorem 6, where
$c_1,c_2,c_3,c_4$ are constants. In particular, this case
induces the Guichard net of Bianchi-type, since we have 
$$\varphi_{xx}=\frac{c_3^2(c_1+c_2)}{2}\sin2\varphi, \hspace{0.5cm} \varphi_{yy}=\frac{c_4^2(c_1+c_2)}{2}\sin2\varphi, \hspace{0.5cm} \varphi_{zz}^c=\frac{c^2(c_1+c_2)}{2}\sin2\varphi$$
and uniqueness of solutions for the evolution equation in $z$ with respect to the initial condition, which we shall study in \S4 (and see Example 5 there). \\

{\sc Proof}. \ \ Let us set \ $(\varphi^1_z)^2=c_1\sin^2\varphi-c_2\cos^2\varphi$. \ Then, we have 
$${\hat A}=-(c_1+c_2)\frac{\varphi_x\cos\varphi}{(\varphi^1_z)^2}, \hspace{1cm} {\hat B}=(c_1+c_2)\frac{\varphi_y\sin\varphi}{(\varphi^1_z)^2}.$$
Then, Theorem 6-(1) is given by \ $\varphi_{xy}(c_1\sin^2\varphi-c_2\cos^2\varphi)=(c_1+c_2)\varphi_x\varphi_y\sin\varphi\cos\varphi$. \ Since \ $(c_1\sin^2\varphi-c_2\cos^2\varphi)'=(c_1+c_2)\varphi'\sin2\varphi$, \ we have 
$$(\varphi_x)^2=\varrho^2(x)(c_1\sin^2\varphi-c_2\cos^2\varphi), \hspace{1cm}(\varphi_y)^2=\sigma^2(y)(c_1\sin^2\varphi-c_2\cos^2\varphi).$$

Let $c_3:=\varrho(x)$ and $c_4:=\sigma(y)$. Then, we have \ $L\varphi=(c_3^2-c_4^2)(c_1+c_2)\sin\varphi\cos\varphi$ \ and 
$$\varphi_x(L\varphi)=\frac{(c_3^2-c_4^2)(c_1+c_2)}{2}\varphi_x\sin2\varphi=\left(\frac{(c_3^2-c_4^2)}{2}(c_1\sin^2\varphi-c_2\cos^2\varphi)\right)_x,$$
$$\varphi_y(L\varphi)=\frac{(c_3^2-c_4^2)(c_1+c_2)}{2}\varphi_y\sin2\varphi=\left(\frac{(c_3^2-c_4^2)}{2}(c_1\sin^2\varphi-c_2\cos^2\varphi)\right)_y.$$
For Theorem 6-(2), we determine
$$S:=\frac{c_3^2-c_4^2}{2}(c_1\sin^2\varphi-c_2\cos^2\varphi+c_2), \hspace{1cm}T:=\frac{c_3^2-c_4^2}{2}(c_1\sin^2\varphi-c_2\cos^2\varphi-c_1). $$
\hspace{\fill}$\Box$\\

{\exa}4. \ \ Let us take \ $(\varphi^1_z)^2(x,y,0)=\zeta(x)\sin^2\varphi(x,y,0)$ \ with any positive function $\zeta(x)$. \ Then, the function $\varphi(x,y,0)$ such that \ $\varphi_x=c_1\sin\varphi$ \ and \ $\varphi_y=c_2\sin\varphi$ \ satisfies the condition (1) and (2) of Theorem 6, where $c_1,c_2$ are constants. In particular, this case induces many metrics ${\hat g}$ determined by any $c_1$, $c_2$ and $\zeta(x)$, of which $(L\psi^c)(x,y,0)$ satisfies Proposition 3.1-(3). \\

{\sc Proof}. \ \ Let us set \ $(\varphi^1_z)^2(x,y,0)=\zeta(x)\sin^2\varphi(x,y,0)$. \ Then, we have 
$${\hat A}=-\frac{1}{2\sin\varphi}\left(\frac{\zeta'}{\zeta}+2\varphi_x\cot\varphi\right), \hspace{1cm} {\hat B}=\frac{\varphi_y}{\sin\varphi}.$$

For Theorem 6-(1), we have \ $\varphi_{xy}=\varphi_x\varphi_y\cot\varphi.$ \ This equation is independent of $\zeta(x)$ and we have \ $\varphi_x=\varrho(x)\sin\varphi$ \ and \ $\varphi_y=\sigma(y)\sin\varphi$. \ 

Now, when we take \ $\zeta(x):=1/4x^2$, \ $\varrho(x):=-1/x$ \ and \ $\sigma(y):=1/y$, \ we obtain $\varphi(x,y,0)$ and the metric ${\hat g}$ of Example 1 in \S2.2. Then, for any $\zeta(x)$, $\varphi(x,y,0)$ satisfies Theorem 6-(2) with respect to $S(x,y)$ and $T(x,y)$ of Example 1, as $\varphi(x,y,0)$ is independent of $\zeta(x)$. Hence, in this case we obtain many examples of ${\hat g}$, of which $(L\psi^c)(x,y,0)$ satisfies Proposition 3.1-(3), by giving arbitrary $\zeta(x)$.    

Here, we assume \ $c_1:=\varrho(x)$ \ and \ $c_2:=\sigma(y)$. \ Then, we have \ 
$$L\varphi=\frac{c_1^2-c_2^2}{2}\sin2\varphi, \ \ \ \varphi_x(L\varphi)=\left(\frac{c_1^2-c_2^2}{2}\sin^2\varphi\right)_x, \ \ \ \varphi_y(L\varphi)=\left(\frac{c_1^2-c_2^2}{2}\sin^2\varphi\right)_y.$$ 
Hence, for Theorem 6-(2), we can take \ 
$$S:=\frac{c_1^2-c_2^2}{2}\sin^2\varphi, \ \ \ \ T:=\frac{c_1^2-c_2^2}{2}(\sin^2\varphi-1)=-\frac{c_1^2-c_2^2}{2}\cos^2\varphi.$$ \hspace{\fill}$\Box$

\section{System of evolution equations and construction of
  Guichard nets}
\label{sec:syst-evol-equat}

In this section, we show that a class of functions $\varphi(x,y,z)$ and $\psi(x,y,z)$ in Theorem 1 is obtained as solutions of a system of evolution equations in $z$ from initial data ${\hat g}$ at $z=0$, which are orthogonal {\it analytic} Riemannian $2$-metrics with constant Gauss curvature $-1$ determined by Theorems 5 and 6.  
Theorems 3, 5 and 6 will be useful to verify this fact.  

Now, we consider the following system of evolution equations in $z$: 
$$\psi_{zz}=(\varphi_{xx}-\varphi_{yy})\sin2\varphi-(\psi_{xx}-\psi_{yy})\cos2\varphi, \eqno{(4.1)}$$
$$\varphi_{zz}=(\varphi_{xx}-\varphi_{yy})\cos2\varphi+(\psi_{xx}-\psi_{yy})\sin2\varphi, \eqno{(4.2)}$$
under a suitable initial condition at $z=0$. 

Now, for the system of (4.1) and (4.2), the initial condition at $z=0$ is obtained from analytic $2$-metrics ${\hat g}$ determined by Theorems 5 and 6: \ Let us choose analytic functions $D(y)$ and $\zeta(x)>0$ in \ $\cos^2\varphi(x,y,0)=1/(1+e^{D(y)})$ \ and \ $(\varphi^c_z)^2(x,y,0)=c^2\zeta(x)\sin^2\varphi(x,y,0)$ \ of Theorem 5 and choose analytic functions $\zeta(x)$, $\eta(y)$ and $\varphi(x,y,0)$ in \ $(\varphi^c_z)^2(x,y,0)=c^2(\zeta\sin^2\varphi-\eta\cos^2\varphi)(x,y,0)$ \ of Theorem 6. Then, an analytic metric 
$${\hat g}={\hat A}^2(x,y)(dx)^2+{\hat B}^2(x,y)(dy)^2 \eqno{(4.3)}$$
is defined from these functions such that ${\hat g}$ is independent of $c$ and has constant Gauss curvature $-1$. Furthermore, $\psi^c(x,y,0)$ and $\psi^c_z(x,y,0)$ are determined for such a metric ${\hat g}$, and all systems of four functions $\varphi(x,y,0)$, $\psi^c(x,y,0)$, $\varphi^c_z(x,y,0)$ and $\psi^c_z(x,y,0)$ depending on $c\neq 0$ satisfy (a), (b) and (c) in Theorem 3 and further satisfy the following (d) and (e):
$$\psi^c_{xy}(x,y,0)=(\varphi_x\varphi_y)(x,y,0), \leqno{(d)}$$
$$[-\Delta\psi^c+(\varphi_x)^2+(\varphi_y)^2+(\varphi^c_z)^2](x,y,0)=
[(\varphi_{xx}-\varphi_{yy})\sin2\varphi-(\psi^c_{xx}-\psi^c_{yy})\cos2\varphi](x,y,0). \leqno{(e)}$$
Conversely, if a metric ${\hat g}$ defines systems of four analytic functions at $z=0$ depending on $c\neq 0$ such that each system satisfies (a), (b), (c), (d) and (e), then ${\hat g}$ is obtained from $\varphi(x,y,0)$ and $\varphi^c_z(x,y,0)$ as above. 
We take systems of four functions determined from such a ${\hat g}$ and $c\neq 0$ as the initial condition for (4.1) and (4.2). \\

{\sc Remark for the Initial Condition}. \ \ Firstly, we note that all 
initial functions at $z=0$ are analytic. This analyticity
for initial functions is necessary because we will apply
the Cauchy--Kovalevskaya Theorem to obtain existence and
uniqueness of solutions of the system (4.1)
and (4.2). 

For $\varphi(x,y,0)$, we can arbitrarily take \ $\varphi(0,0,0)=\lambda$. \ However, when we define \ ${\bar \varphi}(x,y,z):=\varphi(x+a_1,y+a_2,z)$ \ with any constants $a_1 \ {\rm and} \ a_2$, these ${\bar \varphi}$ and $\varphi$ lead to the same Guichard net. 
Hence, we may assume $\varphi(0,0,0)=\pi/4$. 

From $\varphi^c_z(x,y,0)$, we determine $\psi^c_z(x,y,0)$ by (c) as follows:
$$\psi^c_z(x,y,0):=\int^{(x,y,0)}_{(0,0,0)} \{-(\varphi^c_{xz}\cot\varphi)(x,y,0)dx+(\varphi^c_{yz}\tan\varphi)(x,y,0)dy\},$$
that is, $\psi_z^c(0,0,0)=0$ and $\psi^c_z(x,y,0)$ is determined up to the same constant multiple $c$ as $\varphi^c_z(x,y,0)$.

$\psi^c(x,y,0)$ is determined from $(L\psi^c)(x,y,0)$ up to terms of $k(x+y)$ and ${\hat k}(x-y)$ by Theorem 5 and Proposition 3.4. Then, it will be uniquely determined by (d), (e) and the condition that it has no linear term for $x$ and $y$: \ $\psi^c(x,y,0)$ is generally expressed by (d) in the form 
$$\psi^c(x,y,0)=\int_0^x\int_0^y(\varphi_x\varphi_y)(x,y,0)dxdy +X^c(x)+Y^c(y) \eqno{(4.4)}$$
with functions $X^c(x)$, $Y^c(y)$, where we choose $X^c(x)$ and $Y^c(y)$ such that they do not have any linear term for $x$ and $y$.  
Then, $X^c(x)$ and $Y^c(y)$ are uniquely determined by (e), of which fact will be verified in Proposition 4.1 below. 

Thus, we have obtained from an initial data ${\hat g}$ determined by Theorems 5 and 6 a one-parameter family $\{\varphi(x,y,0), \ \psi(x,y,0), \ \varphi^c_z(x,y,0),\psi^c_z(x,y,0)\}$ with parameter $c\neq 0$ as the initial condition. 
Consequently, for a given metric ${\hat g}$, there is a
one-parameter family $\{\varphi^c(x,y,z),\psi^c(x,y,z)\}$ of
solutions for the system of equations (4.1) and (4.2), which
will lead to distinct Guichard nets if $c\neq c'$ (see Theorem 7 below). \\

{\pro}4.1. \ \ {\it Suppose that $\psi^c(x,y,0)$ satisfies (d) and (e). 
Then, $\psi^c(x,y,0)$ is uniquely determined, if it does not have any linear term for $x$ and $y$.}\\

{\sc Proof}. \ In this proof, we omit the $c$ in $\psi^c(x,y,0)$, $X^c(x)$, $Y^c(y)$, etc. 

Now, let $\hat{\psi}(x,y,0)$ be the first integral term in the right hand side of (4.4). Suppose that $\psi(x,y,0)$ has two expressions of \ ${\bar \psi}(x,y,0)=\hat{\psi}(x,y,0)+{\bar X}(x)+{\bar Y}(y)$ \ and \ ${\tilde \psi}(x,y,0)=\hat{\psi}(x,y,0)+{\tilde X}(x)+{\tilde Y}(y)$. \ Then, since 
$${\bar X}''\sin^2\varphi+{\bar Y}''\cos^2\varphi={\tilde X}''\sin^2\varphi+{\tilde Y}''\cos^2\varphi$$
$$=-(1/2)[\Delta\hat{\psi}-\varphi_x^2-\varphi_y^2-\varphi_z^2+(L\varphi)\sin2\varphi-(L\hat{\psi})\cos2\varphi]$$
by (e), we firstly have \ $({\bar X}-{\tilde X})''\sin^2\varphi+({\bar Y}-{\tilde Y})''\cos^2\varphi=0$ \ for $(x,y,0)$.

Next, there are functions $k(x+y)$ and ${\hat k}(x-y)$ such that \ $({\bar \psi}-{\tilde \psi})(x,y,0)=k(x+y)+{\hat k}(x-y)$, \ as $\psi(x,y,0)$ is determined from $(L\psi)(x,y,0)$. Taking derivatives of $({\bar \psi}-{\tilde \psi})(x,y,0)$ by $x$ and $y$ respectively, we have \ $({\bar X}-{\tilde X})''(x)=({\bar Y}-{\tilde Y})''(y)=k''(x+y)+{\hat k}''(x-y)$. \ 

From these two equations, we obtain \ $({\bar X}-{\tilde X})''(x)=({\bar Y}-{\tilde Y})''(y)=0$, \ which shows that $\psi(x,y,0)$ is uniquely determined up to linear terms. 

Finally, we note that, if $\varphi(x,y,0)$ is really a function of two variables $x$ and $y$, then the conclusion of the Proposition follows from only the first equation. \hspace{\fill}$\Box$\\

Now, we define the functions $I^x(x,y,z)$, $I^y(x,y,z)$, $J(x,y,z)$ and $K(x,y,z)$, respectively, by using the solutions $\varphi(x,y,z)$ and $\psi(x,y,z)$ for the system (4.1) and (4.2):
$$ 
\begin{array}{l}
I^x:=\psi_{xz}+\varphi_{xz}\cot\varphi, \hspace{2.5cm} I^y:=\psi_{yz}-\varphi_{yz}\tan\varphi,\\
J:=\psi_{xy}-\varphi_x\varphi_y, \hspace{0.5cm} K:=(L\varphi)\sin2\varphi-(L\psi)\cos2\varphi+\Delta\psi-(\varphi_x)^2-(\varphi_y)^2-(\varphi_z)^2.
\end{array}
\eqno{(4.5)}$$\\

{\pro}4.2. \ \ {We have the following system of equations for any $(x,y,z)$:
\[
\frac{\partial}{\partial z}\left[
\begin{array}{l}
I^x \\
I^y \\
J \\
K 
\end{array}
\right]  
=
\left[
\begin{array}{llll}
\ \ \ \ 0 & \ \ \ \ 0 & -\frac{1}{\sin^2\varphi}\partial/\partial y & \frac{1}{2\sin^2\varphi}\partial/\partial x \\
\ \ \ \ 0 & \ \ \ \ 0 & -\frac{1}{\cos^2\varphi}\partial/\partial x & \frac{1}{2\cos^2\varphi}\partial/\partial y \\
\sin^2\varphi \ \partial/\partial y & \cos^2\varphi \ \partial/\partial x & \ \ \ \ \ 0 & \ \ \ \ \ 0 \\
2\sin^2\varphi \ \partial/\partial x & 2\cos^2\varphi \ \partial/\partial y & \ \ \ \ \ 0 & \ \ \ \ \ 0 \\
\end{array}
\right]
\left[
\begin{array}{l}
I^x \\
I^y \\
J \\
K 
\end{array}
\right]_.
\]}\\

{\sc Proof}. \ We obtain the equations of $(I^x)_z$ and $(I^y)_z$ from the proof of Theorem 4, where we showed that $(I^x)_z=0$ and $(I^y)_z=0$, respectively, are satisfied if and only if the right hand sides of them vanish, by using (4.1) and (4.2).  

For the equation of $J_z$, we firstly define \ ${\hat{I^x}}:=\tan\varphi I^x$ \ and \ ${\hat{I^y}}:=\cot\varphi I^y$. \ We have 
$$({\hat{I^x}})_y+({\hat{I^y}})_x=\displaystyle\frac{1}{\sin\varphi\cos\varphi}[J_z+\varphi_y{\hat{I^x}}-\varphi_x{\hat{I^y}}],$$
then we obtain the equation desired. 

For the equation of $K_z$, we have
$$K_z=[(L\varphi)\sin2\varphi-(L\psi)\cos2\varphi+\Delta\psi-(\varphi_x)^2-(\varphi_y)^2-(\varphi_z)^2]_z $$
$$=(L\varphi)_z\sin2\varphi-(L\psi)_z\cos2\varphi+\Delta\psi_z-2\varphi_x\varphi_{xz}-2\varphi_y\varphi_{yz}$$
$$=(\varphi_{xxz}-\varphi_{yyz})\sin2\varphi+2\psi_{xxz}\sin^2\varphi+2\psi_{yyz}\cos^2\varphi-2\varphi_x\varphi_{xz}-2\varphi_y\varphi_{yz}$$
$$=2\sin^2\varphi(\psi_{xxz}+\varphi_{xxz}\cot\varphi)+2\cos^2\varphi(\psi_{yyz}-\varphi_{yyz}\tan\varphi)-2\varphi_x\varphi_{xz}-2\varphi_y\varphi_{yz}$$
$$=2\sin^2\varphi(I^x)_x+2\cos^2\varphi(I^y)_y.$$ \hspace{\fill}$\Box$\\

The matrix of the right hand side in Proposition 4.2 is a
linear differential operator of first order with respect to
$x$ and $y$, then the system in Proposition 4.2 is regarded
as an evolution equation in $z$. Hence, when we take
solutions $\varphi(x,y,z)$ and $\psi(x,y,z)$ of (4.1) and
(4.2) under the initial condition determined as above, we
obtain \ $I^x\equiv I^y\equiv J\equiv K\equiv 0$ \ for any
$(x,y,z)$ by the uniqueness assertion of the Cauchy--Kovalevskaya, as $I^x(x,y,0)\equiv I^y(x,y,0)\equiv J(x,y,0)\equiv K(x,y,0)\equiv 0$ are satisfied. 

\vspace{2mm}
In the statement and the proof of Theorem 7 below, we assume that $\psi$ does not have any linear term for $x$, $y$, $z$, that is, the initial function $\psi(x,y,0)$ (resp. $\psi_z(x,y,0)$) not only satisfies (d) (resp. (c)) but they, respectively, are also defined by the conditions given in the Remark above.\\

{\thm}7. \ \ {\it Let us take an analytic $2$-metric ${\hat g}$ given in Theorem 5 or obtained by the procedure in Theorem 6.
Let functions $\varphi(x,y,0)$, $\psi^c(x,y,0)$, $\varphi_z^c(x,y,0)$
and $\psi_z^c(x,y,0)$ be a system determined by ${\hat g}$ as above.
We take such a system of functions as the initial condition at $z=0$
for the system (4.1) and (4.2). Then, all solutions $\varphi^c(x,y,z)$
and $\psi^c(x,y,z)$ depending on $c$ satisfy all equations of Theorem
1, that is, each pair $\varphi^c(x,y,z)$ and $\psi^c(x,y,z)$ defines
an evolution of $2$-metrics issuing from ${\hat g}$, which
corresponds to a conformally flat 3-metric with the Guichard condition.  

Conversely, if, for an orthogonal analytic $2$-metric ${\hat g}$ with
constant Gauss curvature $-1$, there is a one-parameter family of
evolutions of $2$-metrics issuing from ${\hat g}$ such that each
evolution corresponds to a conformally flat 3-metric with the Guichard net, then ${\hat g}$ is a metric either in Theorem 5 or obtained by the procedure in Theorem 6.}\\

{\sc Proof}. \ \ 
Let an analytic $2$-metric ${\hat g}$ and a system of
functions $\varphi(x,y,0)$, $\psi(x,y,0)$,
$\varphi_z(x,y,0)$, $\psi_z(x,y,0)$ satisfy the hypotheses
of the theorem. 
Since these four functions given as an initial condition at $z=0$ are analytic, a pair of solutions $\varphi(x,y,z)$ and $\psi(x,y,z)$ for the system (4.1) and (4.2) uniquely exists for each initial condition depending on $c$. 
Hence, we can assume that $\varphi(x,y,z)$ and $\psi(x,y,z)$ satisfy (4.1), (4.2) for any $(x,y,z)$ and also satisfy the initial condition (a), (b), (c), (d) and (e) at $z=0$. 

Then, we obtain  \ $I^x\equiv I^y\equiv J\equiv K\equiv 0$ \ for any $(x,y,z)$ by Proposition 4.2. That is, $\varphi(x,y,z)$ and $\psi(x,y,z)$ not only satisfy (4.1), (4.2) but also satisfy the following equations for any $(x,y,z)$: 
$$\psi_{xz}=-\varphi_{xz}\cot\varphi, \hspace{1cm} \psi_{yz}=\varphi_{yz}\tan\varphi, \eqno{(4.6)}$$
$$\psi_{xy}=\varphi_x\varphi_y,  \eqno{(4.7)}$$
$$\psi_{zz}=(L\varphi)\sin2\varphi-(L\psi)\cos2\varphi=-\Delta\psi+(\varphi_x)^2+(\varphi_y)^2+(\varphi_z)^2.
\eqno{(4.8)}$$

Thus, since the solutions $\varphi(x,y,z)$ and $\psi(x,y,z)$ of the
system (4.1) and (4.2) under our initial condition also satisfy (4.6),
$\varphi(x,y,z)$ and $\psi(x,y,z)$ satisfy all equations in Theorem 1.
In particular, each solution $\{\varphi^c, \psi^c\}$ obtained from
${\hat g}$ and $c\neq 0$ defines an evolution of $2$-metrics issuing
from ${\hat g}$ and the evolution corresponds to a conformally flat 3-metric with the Guichard condition.  

Next, we verify the converse. Let us assume that there is a
one-parameter family of evolutions of $2$-metrics issuing from a
$2$-metric ${\hat g}$ with constant Gauss curvature $-1$ and that each
evolution corresponds to a conformally flat 3-metric $g^c$ with the Guichard net.
Then, ${\hat g}$ determines systems of functions $\varphi(x,y,0)$, $\psi^c(x,y,0)$, $\varphi^c_z(x,y,0)$ and $\psi^c_z(x,y,0)$ depending on $c$ such that each system satisfies (a), (b), (c), (d) and (e) by Theorems 2 and 4. 
On the other hand, by Theorem 3, Corollary 3.3, Theorem 5 and Theorem 6, an orthogonal $2$-metric ${\hat g}$ with constant Gauss curvature $-1$ defines systems of functions $\varphi(x,y,0)$, $\psi^c(x,y,0)$, $\varphi^c_z(x,y,0)$ and $\psi^c_z(x,y,0)$ depending on $c$ such that each system satisfies (a), (b), (c), (d) and (e), if and only if ${\hat g}$ is a metric given in Theorem 5 or obtained by the procedure in Theorem 6. Thus, the converse statement has been proved. 

By these arguments, we have completely verified the Theorem. \hspace{\fill}$\Box$\\

 In general, it seems difficult to solve the system of evolution equations (4.1) and (4.2). 
However, in the case of the Bianchi-type Guichard net, we
may extract the function $\varphi(x,y,z)$ from the initial
metric ${\hat g}$ and then we find $\psi(x,y,z)$ as follows:\\

{\exa}5 \ (Bianchi-type Guichard net). \ \ All functions inducing the Bianchi-type Guichard net are given by $\varphi(x,y,z)=g(ax+by+cz)$, \ $abc\neq 0$, where $g(t)$ is a 1-variable function such that 
$$g''=\alpha\sin 2g, \ \ \ \ (g' )^2=\beta-\alpha\cos 2g \eqno{(4.9)}$$
with constants $\alpha$ and $\beta$. 
Here, in this case,  we study how $\varphi(x,y,z)$ and $\psi(x,y,z)$ are determined from the initial data ${\hat g}$. \\

Firstly, we fix the initial data $\hat{g}$: \ \ let us take 
$$\hat{A}(x,y):=-2a\alpha\frac{\cos g}{g'}(ax+by), \hspace{1cm}
\hat{B}(x,y):=2b\alpha\frac{\sin g}{g'}(ax+by)$$
as in \S2.1. 
Then, the metric \ ${\hat g}=\hat{A}^2(x,y)dx^2+\hat{B}^2(x,y)dy^2$ \ has the constant Gauss curvature $-1$. 
In fact, we can show it from \ $\hat{A}_y/\hat{B}=a(\alpha+\beta)/g'$ \ and \ $\hat{B}_x/\hat{A}=b(\alpha-\beta)/g'$ \ by direct calculation. 

Next, we study the initial condition. 
We have \ $\varphi(x,y,0)=g(ax+by)$ \ from (a). 
The equation (b) implies the following equation: 
$$(\log|\varphi_z|)_x(x,y,0)=(\log |g'|)_x(ax+by), \hspace{1cm}
(\log|\varphi_z|)_y(x,y,0)=(\log |g'|)_y(ax+by).$$
Hence, we have \ $\varphi_z(x,y,0)=cg'(ax+by)$ \ with any constant $c(\neq 0)$. $\psi_z(x,y,0)$ satisfying \ $\psi_z(0,0,0)=0$ \ is determined by (c). 
From \ $\psi_{xy}(x,y,0)=(\varphi_x\varphi_y)(x,y,0)=abg'^2(ax+by)$ \ by (d), we have 
$$\psi(x,y,0)=X(x)+Y(y)+\int_0^{t}ds\int_0^sg'^2(u)du,$$
where $t=ax+by$. 
Then, $X(x)$ and $Y(y)$ are determined by (e): 
$$X(x)=(c_1/2)x^2, \hspace{1.5cm} Y(y)=(c_2/2)y^2, $$
where \ $2c_1=(\alpha+\beta)(-a^2+b^2+c^2)$, \ $2c_2=(\alpha-\beta)(-a^2+b^2-c^2)$. 

Since we have obtained all initial condition $\varphi(x,y,0)$, $\psi(x,y,0)$, $\varphi_z(x,y,0)$, $\psi_z(x,y,0)$ for the system (4.1) and (4.2), a pair of solutions $\varphi(x,y,z)$ and $\psi(x,y,z)$ are uniquely determined. On the other hand, \ $\varphi(x,y,z):=g(ax+by+cz)$ \ satisfies this initial condition and it is known that $g(ax+by+cz)$ induces a conformally flat 3-metric with the Guichard condition. Hence, we may obtain a one-parameter family $\varphi(x,y,z)=g(ax+by+cz)$ with parameter $c(\neq 0)$ from ${\hat g}$, as the partner of $\psi(x,y,z)$.

Now, we shall uniquely determine $\psi(x,y,z)$ from (4.1), (4.6), (4.7) and (4.8). Since \ $\psi_{xz}(x,y,z)=-ac(g''\cot g)(x,y,z)$ \ and \ $\psi_{yz}(x,y,z)=bc(g''\tan g)(x,y,z)$ \ by (4.6), 
we have 
$$\psi_{xz}(x,y,z)=-ac[(\alpha+\beta)-g'^2(ax+by+cz)], \hspace{0.5cm}\psi_{yz}(x,y,z)=bc[(\alpha-\beta)+g'^2(ax+by+cz)].$$ 
Thus, by (4.7), we firstly define ${\hat \psi}$, from which $\psi$ will be produced, by 
$$\hat{\psi}(x,y,z):=X(x)+Y(y)+Z(z)-(\alpha+\beta)(ax)(cz)+(\alpha-\beta)(by)(cz)+\int_0^{t}ds\int_0^sg'^2(u)du,$$
where $t=ax+by+cz$. Then, since \ $\hat{\psi}_{zz}=-\Delta{\hat \psi}+\varphi_x^2+\varphi_y^2+\varphi_z^2$ \ by (4.8) and 
$$\hat{\psi}_{zz}=Z''+c^2g'^2, \ \ \ \ 
-\Delta{\hat \psi}+\varphi_x^2+\varphi_y^2+\varphi_z^2=-X''-Y''+c^2g'^2,$$
we have
$$X''+Y''+Z''=0\Longleftrightarrow X(x)=(c_1/2)x^2, \ \ Y(y)=(c_2/2)y^2, \ \ Z(z)=-[(c_1+c_2)/2]z^2.$$
Since \ $[-\Delta{\hat \psi}+\varphi_x^2+\varphi_y^2+\varphi_z^2](x,y,z)=[L(\varphi)\sin2\varphi-L(\psi)\cos2\varphi](x,y,z)$ \ by (4.1) and (4.8), \ we have
$$-(c_1+c_2)+c^2(\beta-\alpha\cos2g)=\alpha(a^2-b^2)-[c_1-c_2+\beta(a^2-b^2)]\cos2g.$$
Hence, we have again
$$c_1+c_2=-\alpha(a^2-b^2)+\beta c^2, \ \ \ \ c_1-c_2=\alpha c^2-\beta(a^2-b^2).$$

In consequence, we have obtained 
$$\varphi(x,y,z)=g(ax+by+cz),$$ 
$$\psi(x,y,z)=-ac(\alpha+\beta)xz+bc(\alpha-\beta)yz+\frac{c_1}{2}x^2+\frac{c_2}{2}y^2+\frac{c_3}{2}z^2+\int_0^tds\int_0^sg'^2(u)du,$$
where \ $2c_1=(\alpha+\beta)(-a^2+b^2+c^2)$, \ $2c_2=(\alpha-\beta)(-a^2+b^2-c^2)$ \ and $c_1+c_2+c_3=0$. \\

{\small
 
\vspace*{0.5cm}

Francis E.~Burstall\hspace*{5.5cm}       Udo Hertrich-Jeromin\\
Dept of Mathematical Science \hspace*{3cm}Technische Universit\H{a}t Wien, E104\\
University of Bath \hspace*{4.9cm} Wiedner Haupstra{\ss}e 8-10, A-1040 (Austria)\\
Bath BA2 7AY, UK \hspace*{4.65cm} e-mail: uhj@geometrie,tuwien.ac.at \\
e-mail: feb@maths.bath.ac.uk

\vspace{2mm}
Yoshihiko Suyama\\
Dept of Applied Mathematics\\
Fukuoka University, Fukuoka 814-0180 (Japan)\\
e-mail: suyama@fukuoka-u.ac.jp  }

\end{document}